\documentclass[12pt]{amsart}
\usepackage{enumerate,amssymb,amsfonts,latexsym,psfrag}
\usepackage{amscd,amsmath}
\usepackage{bm}
\usepackage[height=9in,width=6.5in,centering]{geometry}
\usepackage[dvips]{graphicx}
\usepackage{overpic}
\usepackage[usenames,dvipsnames]{color}
\usepackage{datetime}
\newcommand{\now}{\today\ at \currenttime} 
\usepackage[implicit=true,
   colorlinks=true,linkcolor=blue,citecolor=PineGreen,urlcolor=red]{hyperref}

\newif\ifWantNotes\WantNotestrue
\newif\ifWantOldNotes\WantOldNotesfalse
\WantNotesfalse\WantOldNotesfalse
\providecommand\NoteGuts[1]{%
\fbox{\parbox{\marginparwidth}{\baselineskip=6pt\par\noindent\raggedright\sffamily{\footnotesize #1}}}}
\providecommand\note[1]{%
\ifWantNotes\normalmarginpar\marginpar{\NoteGuts{#1}}\fi}
\providecommand\rnote[1]{%
\ifWantNotes\reversemarginpar\marginpar{\NoteGuts{\color{blue}#1}}\fi}
\providecommand{\noteC}[1]{%
\ifWantNotes\leavevmode\vadjust{\smash{\llap{\NoteGuts{#1}\kern\marginparsep}}}\fi}

\ifWantOldNotes\WantNotestrue\fi
\providecommand{\Oldnote}[1]{\ifWantOldNotes\note{\tiny(old) #1}\else\relax\fi}
\providecommand{\OldnoteC}[1]{\ifWantOldNotes\noteC{\tiny(old) #1}\else\relax\fi}
\providecommand{\Oldrnote}[1]{\ifWantOldNotes\rnote{\tiny(old) #1}\else\relax\fi}

\newcommand{\C}{{\mathbb C}}
\newcommand{\D}{{\mathbb D}}
\newcommand{\R}{{\mathbb R}}
\newcommand{\CC}{{\mathcal C}}
\newcommand{\HH}{{\mathcal H}}
\newcommand{\OO}{{\mathcal O}}

\newcommand{\RR}{{\mathcal R}}
\renewcommand{\SS}{{\mathcal S}}
\newcommand{\VV}{{\mathcal V}}
\newcommand{\NN}{{\mathcal N}}

\newcommand{\Bas}{\operatorname{Basin}}
\newcommand{\Arg}{\operatorname{Arg}}
\newcommand{\card}{\operatorname{card}}

\newcommand{\length}{\operatorname{length}}
\newcommand{\inter}{\operatorname{int}}
\newcommand{\Ree}{\operatorname{Re}}
\newcommand{\Imm}{\operatorname{Im}}

\usepackage{mathptmx}
\usepackage[arrow, matrix,curve]{xy}
\newcommand{\dfn}[1]{{\textbf{\emph{#1}}}}
\newcommand{\vargamma}{{\boldsymbol{\gamma}}}
\newcommand{\ThePath}{\vargamma}
\newcommand{\Roots}{\RR}
\newcommand{\Vor}[1]{\operatorname{Vor}{\!\bigl(#1\bigr)}}
\newcommand{\vor}[1]{\operatorname{vor}{\!\bigl(#1\bigr)}}
\newcommand{\Set}[1]{\left\{{#1}\right\}}
\newcommand{\st}{\mid}
\newcommand{\ST}[1]{\left| {#1} \right.}
\newcommand{\ray}[1]{\ell_{#1}}
\newcommand{\fhat}{\widehat{f}}
\newcommand{\hatray}[1]{\widehat{\ell}_{#1}}
\newcommand{\slit}[1]{\sigma_{#1}}
\newcommand{\BigO}[1]{\OO\kern-.2em\left({#1}\right)}
\newcommand{\PDone}{{\mathcal{P}\kern-0.25em}_{d,1}}
\newcommand{\PDoneBar}{\overline{\PDone}}
\newcommand{\Nsteps}[1]{\mbox{\large\#}_{#1}}
\newcommand{\Avgsteps}[1]{\overline{\Nsteps{#1}}}
\newcommand{\Lambar}{\overline{\Lambda}}
\newcommand{\Cost}[1]{\operatorname{Cost}(#1)}
\newcommand{\VisG}{{\mathcal{G}}}
\newcommand{\Good}[1]{\mathrm{Good}_{#1}}
\newcommand{\Bad}[1]{\mathrm{Bad}_{#1}}
\newcommand{\mystrut}{\relax}
\newcommand{\LiftSeg}[1]{\bigl[\kern-.55em\bigl[%
     {\mystrut{\,#1\,}}\bigr]\kern-.55em\bigr]}
\newcommand{\Infl}{\mathcal{I}}
\newcommand{\dist}[1]{\left\|{#1}\right\|}
\newcommand{\mycaption}[1]{\caption{\small #1}}
\newcommand{\Th}{\ensuremath{^{\text{th}}}}
\newcommand{\EqRef}[1]{eqn.~\kern-.08em(\ref{#1})}
\newcommand{\alphaZero}{\ensuremath{3-\sqrt{8}}}

\newtheorem{thm}{Theorem}[section]
\newtheorem*{main-thm}{Main Theorem}
\newtheorem{cor}[thm]{Corollary}
\newtheorem{lem}[thm]{Lemma}
\newtheorem{prop}[thm]{Proposition}
\newtheorem{defn}[thm]{Definition}

\newtheoremstyle{myRemark}
    {1.5ex plus 4pt minus 2pt}   
    {1.5ex plus 4pt minus 2pt}   
    {\normalfont}                
    {}                           
    {\itshape}                  
    {.}                          
    {.5em}                       
    {}  

\theoremstyle{myRemark}
\newtheorem{rem}[thm]{Remark}

\newtheorem{problem}[thm]{Question}

\numberwithin{equation}{section}
\numberwithin{figure}{section}

\makeatletter
\let\oldfigure\figure
\def\figure{\@ifnextchar[\figure@i \figure@ii}
\def\figure@i[#1]{\addtocounter{thm}{1}\oldfigure[#1]}
\def\figure@ii{\addtocounter{thm}{1}\oldfigure}
\let\oldtable\table
\def\table{\@ifnextchar[\table@i \table@ii}
\def\table@i[#1]{\addtocounter{thm}{1}\oldtable[#1]}
\def\table@ii{\addtocounter{thm}{1}\oldtable}
\makeatother

\DeclareRobustCommand{\CIRC}{\raisebox{-.1ex}{\includegraphics[height=1.2ex]{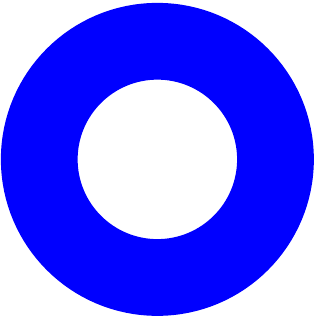}}}
\DeclareRobustCommand{\LINE}{\raisebox{.3ex}{\includegraphics[width=.8em]{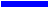}}}
\DeclareRobustCommand{\CROSS}{\raisebox{-.1ex}{\includegraphics[height=1.2ex]{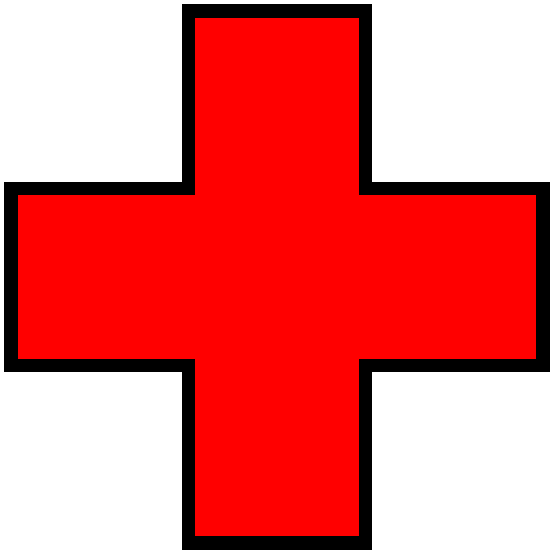}}}

\begin{document}

\title[Geometry of Polynomials and Root-Finding \ifWantNotes(\now)\fi]
   {Geometry of Polynomials and Root-Finding via Path-Lifting}
\author{Myong-Hi Kim}
\address{Myong-Hi Kim \\
Mathematics, Computers \& Information Science \\
SUNY at Old Westbury\\
Old Westbury, NY 11568 USA}
\email{kimm@oldwestbury.edu}
\thanks{{\it Acknowledgements}: Part of this work was done while Myong-Hi
  Kim was visiting Stony Brook University; we are grateful for the
  support and hospitality of the Mathematics Department and the Institute for
  Mathematical Sciences. \\
Marco Martens was supported in part by the National Science Foundation under
the grant DMS-1600554.\\
  A previous version of this article was circulated as ``Bounds for the Cost of
Root Finding.''
  The current version will appear in {\sl Nonlinearity} in modified form.
} 
\author{Marco Martens}
\address{Marco Martens\\
Institute for Mathematical Sciences\\
Stony Brook University\\
Stony Brook, New York 11794 USA}
\email{marco@math.stonybrook.edu}
\author{Scott Sutherland}
\address{Scott Sutherland\\
Institute for Mathematical Sciences\\
Stony Brook University\\
Stony Brook, New York 11794 USA}
\email{scott@math.stonybrook.edu}


\keywords{Root-finding, alpha~theory, Newton's~method, Voronoi~region, 
path-lifting, branched covering, approximate zero, homotopy method}

\subjclass{Primary 65H05; Secondary 30C15, 37F10, 52C20, 57M12, 68Q25}


\begin{abstract} 
Using the interplay between topological, combinatorial, and geometric
properties of polynomials and analytic results (primarily the covering
structure and distortion estimates), 
we analyze a path-lifting method for finding approximate zeros,  
similar to those studied by Smale, Shub, Kim, and others.  
Given any polynomial, this simple algorithm always converges to a root,
except on a finite set of initial points lying on a circle of a given
radius. 

Specifically, the algorithm we analyze consists of iterating
$$z - \frac{f(z)-t_kf(z_0)}{f'(z)}$$ 
where the $t_k$ form a decreasing sequence of real numbers and $z_0$ is chosen
on a circle containing all the roots.
We show that the number of iterates required to locate an
approximate zero of a polynomial $f$ depends only on $\log|f(z_0)/\rho_\zeta|$
(where $\rho_\zeta$ is the radius of convergence of the
branch of $f^{-1}$ taking $0$ to a root $\zeta$) 
and the logarithm of the angle between $f(z_0)$ and certain
critical values.  
\OldnoteC{changed ``inverse of the angle'' to ``reciprocals of these angles'';
  reordered last two sentences.}   
Previous complexity results for related algorithms depend
linearly on the reciprocals of these angles. 
Note that the complexity of the algorithm does not depend directly on the
degree of $f$, but only on the geometry of the critical values.

Furthermore,
for any polynomial $f$ with distinct roots, the average number of steps 
required over all starting points taken on a circle containing all  the
roots is bounded by a constant times the average of $\log(1/\rho_\zeta)$.
The average of $\log(1/\rho_\zeta)$ over all polynomials $f$ with $d$ roots in
the unit disk is $\BigO{d}$.
\OldnoteC{changed ``average number'' to ``expected number''}
This algorithm readily generalizes to finding all roots of a
polynomial (without deflation); doing so increases the complexity by a
factor of at most $d$.
\end{abstract} 

\maketitle


\thispagestyle{empty}

\setcounter{tocdepth}{1}

\vspace{-.5in}   
\renewcommand{\contentsnamefont}{\small\scshape}
\tableofcontents
\addtocontents{toc}{\protect\small} 

\section{Introduction}\label{intro}

We analyze a path-lifting method called
the \dfn{$\alpha$-step method} 
(see page~\pageref{plm} in Section~\ref{alg} for specifics),
which locates an approximate zero (see Definition~\ref{approxzero}) 
for a 
complex polynomial $f(z)$;
from an approximate zero, Newton's method converges quadratically to a
root.
For any polynomial, the $\alpha$-step method converges everywhere except on a
finite set of starting points lying on a circle of given radius. 
This is established in this paper, but also follows from \cite[Thm~5A,5B]{K2}. 

We consider monic polynomials of degree $d$ with distinct roots in the unit
disk, and denote the set of all such polynomials by $\PDone$.
Our main results bound the number of iterations required to locate an
approximate zero in three contexts:
 we bound the number of steps needed to locate an approximate zero
 starting from any point $z_0$ on a circle containing all the roots; 
 we compute the average number of steps over the circle of initial
 points;
 we average this quantity over all polynomials in $\PDone$ to get 
 a bound in terms of the degree.
 These bounds apply to all roots of a given polynomial, and can be
 applied to locate all of the roots with a $d$-fold increase in effort. 

\smallskip
While we analyze the complexity of the $\alpha$-step method, it is not our
primary goal to demonstrate that this achieves the optimal bound.  Indeed,
there are certainly other algorithms with a lower worst-case arithmetic 
complexity (at least for finding $\epsilon$-roots) such as that of Pan 
\cite{Pan2} which achieves the nearly optimal bound, or of Renegar
\cite{Renegar} or Kim-Sutherland \cite{KS}.
Some further remarks discussing the arithmetic complexity of
these and other related methods appear toward the end of this section.



Rather, our goal is to examine how the underlying geometry of a polynomial
can be exploited in root-finding methods. 
Tight upper and lower bounds on the radius of convergence of the inverse of
an analytic map are given by $\alpha$-theory; these are useful in
understanding the geometry of the polynomial.
Since the $\alpha$-theory also applies in the multivariate case, it is our
belief that a better understanding of the univariate case  will be
aid in understanding the case of several variables.

\subsection*{Background}\label{intro_background}
We now discuss some background related to path-lifting methods in general.

Path-lifting methods are a class of homotopy methods, and are also
refered to as ``modified Newton's method'' or ``guided Newton's
method''. 
In such methods, it is often useful to distinguish between the domain and
range, so we have 
$$f:\C_{source} \rightarrow \C_{target};$$
the goal is to lift a path $\ThePath$ lying in $\C_{target}$ to one in
$\C_{source}$ leading from an initial point $z_0$ to a root $\zeta$.
Numerically, this is accomplished by constructing a sequence of points
$z_j\in\C_{source}$ via analytic continuation, in such a way that each
$f(z_j)$ approximates the path $\ThePath$ in $\C_{target}$ and gives an
approximation of the lift $f^{-1}(\ThePath)$ in $\C_{source}$.
%
In  this form, such methods were introduced 
by Shub and Smale (see, for example \cite{SS_Complexity2} or \cite{Sm85}), 
although one could argue (as Smale points out in \cite{Sm81}) that 
in some sense this idea goes back to Gauss.   See
\cite{Renegar} and the references therein, as well as \cite{KS}.
The series \cite{Bezout1, Bezout2, Bezout3, Bezout4, Bezout5, Bezout6, 
Bezout7} discusses related methods for systems of polynomial
equations, as does \cite{BeltranPardo}.
A survey of complexity results for solving polynomial equations in one
variable can be found in \cite{Pan}; see also \cite{batra_survey}.
\medskip

The difficulty of computing a local branch of $f^{-1}$ along a path $\ThePath$
in the target space is related to how close $\ThePath$ comes to a critical
value of $f$. 
However, not all critical values of $f$ are relevant: 
if we fix a branch of $f^{-1}$, then for points $y\in\ThePath$ the only
critical points that have an impact are those $c$ for which $f(c)$
lies on the boundary of the largest disk where $f^{-1}(y)$ is analytic.
Consequently, it is useful to
factor $f$ through the (branched) Riemann surface $\SS$ for $f^{-1}$,
giving
$$\xymatrix{
  \C_{source} \ar[r]^-{\fhat} \ar@/_{1pc}/[rr]_f 
& \SS       \ar[r]^-{\pi}
& \C_{target}
} .$$

\begin{figure}[tbp]
\begin{center}
  \begin{overpic}[height=.3\hsize,tics=5]{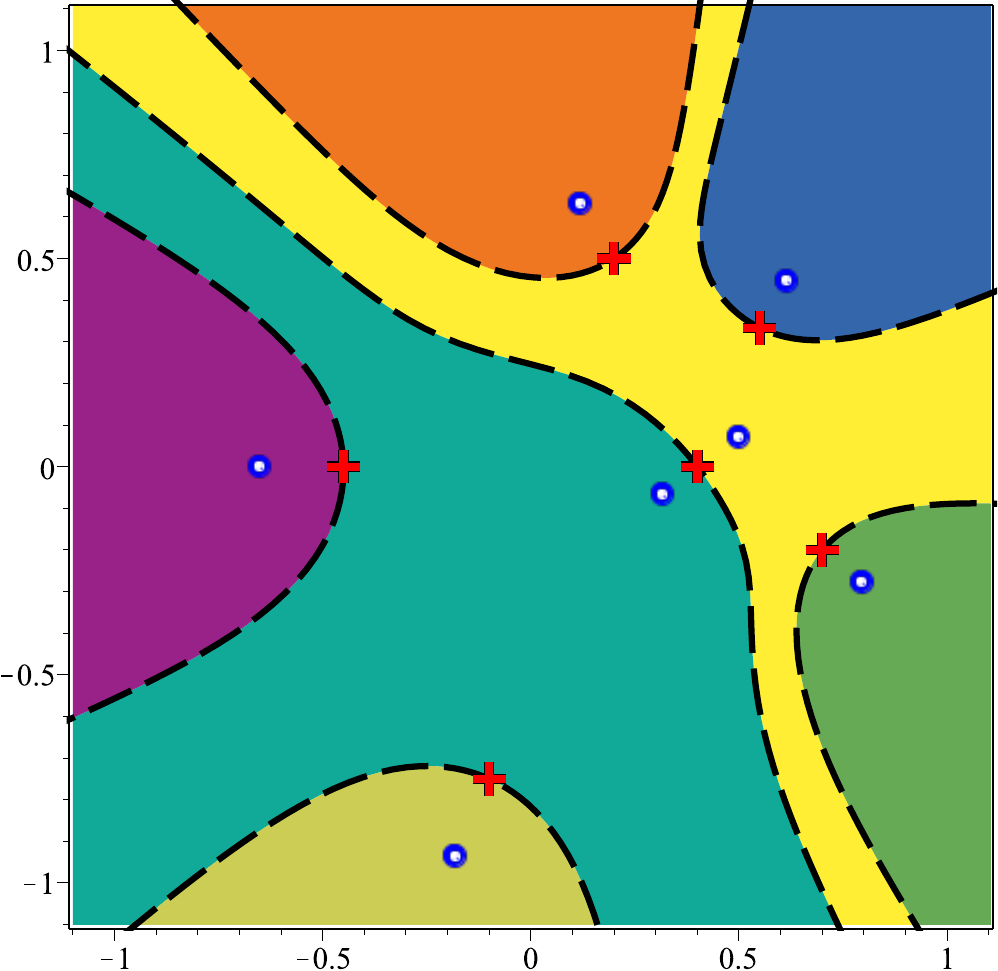} 
    \put(100,55){$B_1$}
    \put(100,90){$B_2$}
    \put(40,100){$B_3$}
    \put(100,30){$B_4$}
    \put(63,-3){$B_7$}
    \put(40,-3){$B_6$}
    \put(-3,40){$B_5$}
  \end{overpic}\quad
  {\vbox to .3\hsize{\vfil\hbox to .1\hsize{\hfil
     $\displaystyle\mathrel{\mathop{\longrightarrow}^{\fhat}}$
     \hfil}\vfil}}
  \begin{overpic}[height=.3\hsize,tics=5]{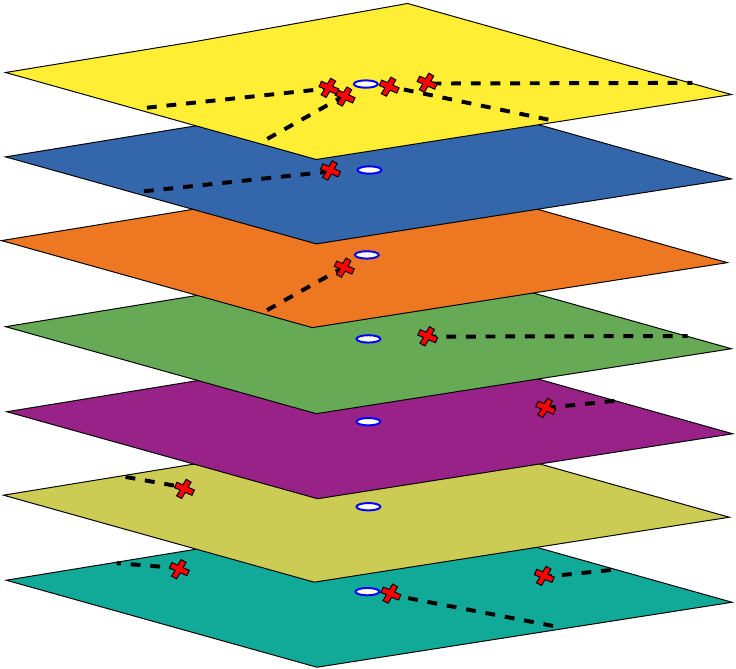} 
    \put(100,77){$\fhat(B_1)$}
    \put(100,65){$\fhat(B_2)$}
    \put(100,52){$\fhat(B_3)$}
    \put(100,41){$\fhat(B_4)$}
    \put(100,29){$\fhat(B_5)$}
    \put(100,17){$\fhat(B_6)$}
    \put(100,5) {$\fhat(B_7)$}
  \end{overpic}\qquad
\end{center}

\OldnoteC{figure \ref{RSpicture} is new}
\mycaption{\label{RSpicture}
For a degree~7 polynomial $f$, on the right is a depiction of
the branched Riemann surface $\SS$ as stack of 7 slit planes.
One side of each slit (indicated by a dashed line)  is joined to the other side
of the parallel slit in a plane above or below it, and
vice-versa. Each slit joins a branch point (indicated by a cross \CROSS) to
infinity. 
On the left, $\C_{source}$ is shown, colored by the corresponding
region of $\SS$; the map~$\fhat$ sends $\C_{source}$ to $\SS$. 
Each critical point of $f$ is marked by a cross, and the preimages of the
slits which terminate at each critical point are indicated by dashed lines.
For reference, the roots of $f$ and their images under $\fhat$
are indicated by circles (\CIRC).  The projection map 
$\pi:\SS\rightarrow\C_{target}$ identifies 
a point in one of the sheets of $\SS$ with all other points directly above and
below it; $\C_{target}$ is not shown in this figure.
} 
\end{figure}

\goodbreak
Denoting the set of critical points $c_j$ of $f$ by $\CC_f$ and the branch
points of $\SS$ by $\VV_f$, we require the map~$\fhat$ to be a biholomorphism
from $\C\smallsetminus\CC_f$ to $\SS\smallsetminus\VV_f$ and a bijection
from $\CC_f$ to $\VV_f$.  Furthermore, the projection~$\pi$ is a $d$-fold
branched cover, and we can choose a metric on $\SS$ so that $\pi$ is a local
isometry away from the branch points. See Figure~\ref{RSpicture}.

The construction of the branched Riemann surface $\SS$ for~$f^{-1}$ is quite
standard, going back to Riemann's dissertation \cite{RiemannThesis},
although often it is presented somewhat abstractly.  Many readers will be 
familiar with the corresponding surfaces for the logarithm and square root;
the explicit view taken here of $\SS$ as a collection of copies of $\C$
identified along slits is similar to the one in \cite[\S10.4]{GreeneKrantz}
or \cite[\S6.1]{Marsden}, to which we refer the interested reader.
Note that each point of $\SS$ corresponds to a pair $(z,w)$ with
$z\in \C_{source}$ and $w\in \C_{target}$, and $w=f(z)$.
It is often helpful to think of the path $\ThePath$ as lying in~$\SS$
rather than in $\C_{target}$; this is possible since for any ray which
avoids $\VV_f$ there is a neighborhood~$U$ containing it which is isometric
to its projection  
$\pi(U)$ in~$\C_{target}$.

\begin{figure}[htpb]
\begin{center}
\begin{overpic}[width=.5\hsize,tics=5]{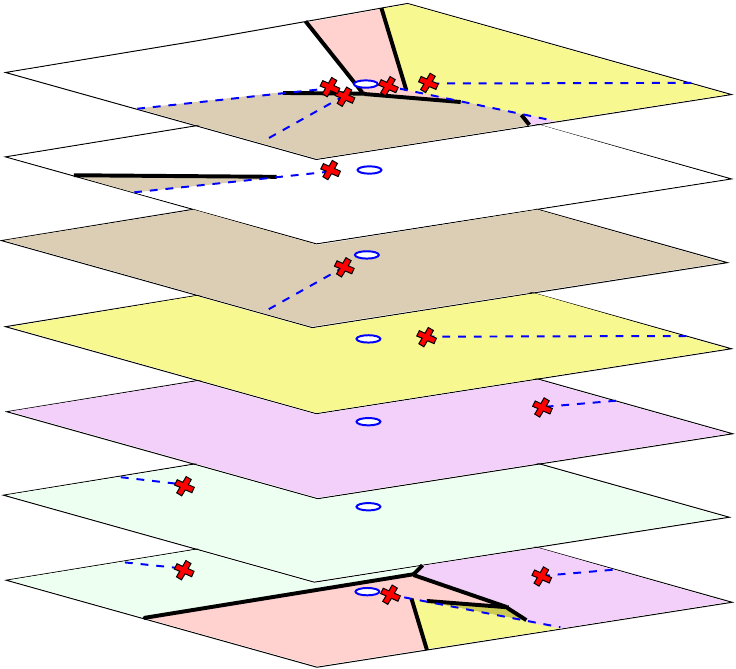}
\put(96,79.5){\tiny $\slit{6}$}
\put(96,45)  {\tiny $\slit{6}$}
\put(13,26.8)  {\tiny $\slit{5}$}
\put(13,15.5)  {\tiny $\slit{5}$}
\put(85,37)  {\tiny $\slit{4}$}
\put(85,14)  {\tiny $\slit{4}$}
\put(40,72)  {\tiny $\slit{3}$}
\put(40,49)  {\tiny $\slit{3}$}
\put(15,74.5)  {\tiny $\slit{2}$}
\put(15,63)  {\tiny $\slit{2}$}
\put(77,3)   {\tiny $\slit{1}$}
\put(65,75){\tiny $\xymatrix @C=15ex @R=1pt{\relax&\slit{1}\ar@/_{.5pc}/[l]}$}
\put(22,89){\tiny $\xymatrix @C=8ex @R=1pt{\Vor{v_1}\ar@/^{.5pc}/[dr]&\relax\\&\relax}$}
\put(40,5){\tiny $\Vor{v_1}$}
\put(16,80){\tiny $\Vor{v_2}$}
\put(40,62){\tiny $\Vor{v_2}$}
\put(49,73){\tiny $\Vor{v_3}$}
\put(16,55){\tiny $\Vor{v_3}$}
\put(-7,65){\tiny $\xymatrix @C=5ex{\Vor{v_3}\ar@/_{.2pc}/[r]&\relax}$}
\put(69,6){\tiny $\xymatrix @C=15ex @R=1pt{\relax&\Vor{v_3}\ar@/_{.5pc}/[l]}$}
\put(16,33){\tiny $\Vor{v_4}$}
\put(62,12){\tiny $\Vor{v_4}$}
\put(70,74){\tiny $\xymatrix @C=18ex @R=1pt{%
           \\ \relax&\Vor{v_4}\ar@/^{.7pc}/[ul]}$}
\put(13,10){\tiny $\Vor{v_5}$}
\put(40,17){\tiny $\Vor{v_5}$}
\put(16,45){\tiny $\Vor{v_6}$}
\put(60,83){\tiny $\Vor{v_6}$}
\put(60,2){\tiny $\Vor{v_6}$}

\end{overpic}
\OldnoteC{Adjusted figure and caption}
\mycaption{\label{vorsheets}
As in Fig.~\ref{RSpicture}, the surface $\SS$ for a degree~7 polynomial is
shown as a stack of seven slit planes, but with Voronoi domains shaded.  
Each sheet is $\fhat(\Bas(\zeta_i))$ for the root $\zeta_i$, 
and is slit along $\slit{v_j}$ (dashed lines), which terminate at the branch
points $v_j\in \VV_f$ (indicated by crosses \CROSS).  
The circles (\CIRC) in each sheet indicate $\pi^{-1}(0)$.
For readability, $\slit{v_j}$ is labeled as $\slit{j}$ in the figure.  
The Voronoi domains of each of the $v_j$ are the labeled regions in the same
shade, with boundaries marked by heavy solid lines (these regions will pass
through slits $\slit{v_k}$ and appear in two or more sheets).
Note that while $\Vor{v_j}$ may enter many sheets, the projection is at most
2-to-1, as in Cor.~\ref{twotoone}. See also Figure~\ref{voronoipic}.
}
\end{center}
\end{figure}

In order to explicitly describe which critical values are relevant for the
path-lifing process, it
is helpful to introduce the Voronoi decomposition of $\SS$ relative to
the branch points $\VV_f$. 
That is, for  each branch point $v$ of~$\SS$, the Voronoi domain
$\Vor{v}$ is the set of 
points in $\SS$ which are closer to $v$ than any other branch point of
$\SS$.
See Figure~\ref{vorsheets}. 
Note that $y\in\Vor{v}$ exactly when $\|v-y\|$ is the radius of convergence
of $\fhat^{-1}$ at $y$.
We show in \S\ref{branched} that the projection map $\pi$ restricted
to any single $\Vor{v}$ is at most $(m+1)$-to-one, where $m$ is the
multiplicity of the critical point of $f$ corresponding to $v$ (hence
the projection $\pi$ is generically at most 2-to-one on $\Vor{v}$).  
When lifting a path $\ThePath$, the number of steps required 
depends directly on the size of a neighborhood of $\ThePath$ on which
a branch of $f^{-1}$ can be defined.  If we think of $\ThePath$ as
lying in $\SS$, then the size of this neighborhood is the distance between
$\ThePath$ and branch points $v_j$ for which $\ThePath$ intersects $\Vor{v_j}$.
We refer to such a critical value $f(c_j)=\pi(v_j)\in\C_{target}$ as
\dfn{relevant} or say that it \dfn{influences} the points on $\ThePath$.

\medskip
As noted earlier, in a path-lifting method we choose a path $\ThePath$
in the target space which connects~$0$ to a point $w_0$ for which we
know a point $z_0$ with $f(z_0)=w_0$.  Path-lifting methods attempt to
calculate a sequence of points $\{z_j\}$ so that $f(z_j) \approx w_j$,
and terminate when a point $z_n$ is sufficiently close to a root
of~$f$. 

Typically the chosen path $\ThePath\in\C_{target}$ is a segment of a ray, and 
we use such paths here. 
It is common (e.g.\ \cite{SS_Complexity2}, \cite{KS}) to choose the
guide points $w_j$ to be of the form $h^j w_0$ for some $h<1$, and
then use one step of Newton's method to obtain $z_{j+1}$ from $z_j$
as
  $z_{j+1} = z_j - (f(z_j)-w_{j+1})/f'(z_j)$.
To ensure convergence, one must choose the $w_j$ sufficiently closely
spaced along $\ThePath$; exactly how close depends strongly on the
size of a neighborhood around~$\ThePath$ on which a branch of $f^{-1}$
can be defined via analytic continuation.

\smallskip
While the ultimate goal of root-finding is typically to find 
a point that lies within an $\epsilon$-ball of some root $\zeta$ of $f$
(called an \dfn{$\epsilon$-root} of $f$), we instead focus on the
problem of locating an \dfn{approximate zero} of $f$.  This notion is
was introduced by Smale (see \cite{Sm81}): a point $z^*$ is an
approximate zero for $f$ if Newton iteration converges at a definite,
rapid rate to a root of $f$ when begun at $z^*$. (See
Definition~\ref{approxzero} for a precise statement.)
   From an approximate zero an $\epsilon$-root for any desired value of
$\epsilon$ can be produced rapidly, with $\BigO{\log|\log\epsilon|}$
iterations of Newton's method (see \cite{Sm85}).

Unlike $\epsilon$-roots, the set of approximate zeros is an intrinsic feature
of a polynomial and does not depend on an externally imposed quantity
$\epsilon$.  We restrict our attention to polynomials with
distinct roots, so approximate zeros always exist for each root $\zeta$.
See also Remark~\ref{epsilonRootFinding} concerning locating
$\epsilon$-roots.

\medskip
Rather than using a regular spacing for the target points $w_j$ in
the path-lifting process, the \mbox{$\alpha$-step} method considered
here selects the points $w_j$ adaptively, 
spacing them as far apart as possible while ensuring 
that at each step $z_j$ is an approximate zero for the function
$f(z)-w_{j+1}$ (and hence $z_{j+1}$ is a good approximation for $w_{j+1}$
with known error bounds). The algorithm terminates when~$z_n$ is an 
approximate zero for $f(z)$.
The tool we use to detect approximate zeros is the Kim-Smale
$\alpha$~function: if $\alpha_f(z) < \alphaZero$, then $z$ 
is an approximate zero for $f$.  See the beginning of
Section~\ref{alg} for further details regarding the $\alpha$~function
and approximate zeros, as well as the specifics of the $\alpha$-step method.
  
\subsection*{Main Results}\label{intro_results}
Our first main result gives an upper bound on the number of steps required
by the $\alpha$-step algorithm to converge to an approximate zero of some
root~$\zeta$ of~$f$, starting from an initial point $z_0 \in \Bas(\zeta)$.
The set $\Bas(\zeta)$ is the collection of all points which converge to the 
root~$\zeta$ under the Newton flow (see Section~\ref{NewtFlow}).
The union of these basins over all roots has full measure; in fact, the
complement is a collection of $d-1$ curves.

The upper bound in the theorem depends on several
quantities closely related to the geometry of the critical values of
$f$.  Specifically, the number of steps required depends on the radius of
convergence $\rho_\zeta$ of the branch of $f^{-1}$ taking $0$ to
$\zeta$ (that is, the norm of some critical value $|f(c_\zeta)|$; this  is
closely related to $f'(\zeta)$), as well as on the angle that the path 
$\ThePath$ makes with the relevant critical values $f(c_j)$ 
(these angles are denoted $\theta_j$ in the statement below) and on the
length of the path (which is $|f(z_0)|$).   As noted earlier, a
critical value $f(c_j)$ is relevant if the corresponding lift of the path
$\ThePath$ to $\SS$ intersects the Voronoi domain of $\fhat(c_j)$.
The appearance of $\rho_\zeta$ in our estimates is not surprising, since the
radius of a disk of approximate zeros about a root $\zeta$ is at most
$\rho_\zeta$. 

Note that the number of steps will be
infinite if either the root $\zeta$ is a multiple root (in which case
$\rho_\zeta = 0$) or there is a relevant critical value $f(c_j)$ lying on
the path $\ThePath$ (in which case $\theta_j=0$).  Since we are
working in $\PDone$, the roots are all distinct (so $\rho_\zeta >0$)
and there are at most $d-1$ paths $\ThePath$ which can contain critical values.

Precise definitions of the terms in the theorem below will take
some time to set up, but we hope the informal discussion above will
give the reader a sense of their meaning. 

\newcounter{thmcount}
\theoremstyle{plain}
\newtheorem*{thmA}{Theorem 1}
\newtheorem*{thmB}{Theorem 2}
\newtheorem*{thmC}{Theorem 3}
\newtheorem*{thmD}{Theorem 4}

\begin{thmA}
Let $f\in\PDone$, and let $z_0$ be an initial point for the
$\alpha$-step path lifting algorithm with $|z_0| > 1$.  Then the
maximum number of steps required for the algorithm to produce an
approximate zero in $\Bas(\zeta)$ starting from $z_0$ is 
$$
\Nsteps{f}(z_0) \le 67 \cdot
 \left( \log\frac{|f(z_0)|}{\rho_\zeta} +\log{40}
      + \sum_{j=1}^{\beta^+(z_0)} 
        ( 3-2\log|\theta_j| )
 \right).
$$
\end{thmA} 

Observe that Theorem~\ref{PointwiseCostBound} implies that for
$f\in\PDone$, the $\alpha$-step algorithm 
converges to a root $\zeta$ for every initial point $z_0$ as long as
$\theta_j\ne 0$. Thus, the algorithm can only fail
for at most $2d-2$ initial points $z_0$ on a circle of fixed radius
larger than 1. 
See also Remark~\ref{algConverges}. 

The details of this theorem are established in Section~\ref{cost}. It is
worth noting that for every polynomial, 
the expected number of relevant critical values 
($\beta^+(z_0)$)
is no more than~2 (as shown in Proposition~\ref{ints}); 
a relation between $\rho_\zeta$ and $f'(\zeta)$ is given in
Lemma~\ref{gammaVsRho}.
 
\medskip
We should emphasize that in the literature 
the dependence on the reciprocal of the angle $|\theta_j|$ is linear (see
\cite{SmaleComplexityThy} for an overview), while in
Thm.~\ref{PointwiseCostBound} the dependence is logarithmic.
Beltr\'an and Shub have recently shown (see~\S7 of \cite{BeltranShub13} or
\cite{BeltranShub10}) the existence of homotopy methods whose number of steps
depends logarithmically on a quantity comparable to our $\theta_j$ (in
projective space), but currently  there is no known constructive method to
produce the necessary path.
Since our paths are line segments in the target space,
this is a significant improvement.

\bigskip
For any fixed polynomial $f$, our second main result gives a bound on
the expected value of the number of steps required when an initial
point is taken on the circle of radius $1+1/d$ (with uniform measure on
the circle).  
This is established in Section~\ref{average}. 

\begin{thmB}
Let $f:\C\to \C$ be a monic polynomial with distinct roots $\zeta_i$ in the
unit disk. 
Let $\Avgsteps{f}$ be 
the average number of steps
required by the $\alpha$-step algorithm to locate an approximate zero for
$f$, where the average is taken over starting points on the circle of radius
$1+1/d$ with uniform measure.
Then
$$
\Avgsteps{f} \le 
134 \left(\frac{1}{d}\sum_{i=1}^{d} \log\frac{1}{\rho_{\zeta_i}} + 6.2\right)
.$$

\end{thmB}

We wish to emphasize that for a specific polynomial $f$, this bound
does not depend directly on the degree, 
but only on the arrangement of the critical values 
(or, more precisely, on the geometry of the branched surface $\SS$).
While $\log 1/\rho_\zeta$ is not bounded above or below for $f\in\PDone$,
its average value grows no more than linearly in the degree of $f$ 
(as stated in Theorem~\ref{LambdaLinearGrowth}, established in
Section~\ref{CostVsDegree}).  

As is apparent in Theorem~\ref{AverageCostBound}, the sum of the
logarithms of the $\rho_\zeta$  plays a crucial role in the
estimates.  Indeed, this quantity is a direct measurement of the
difficulty of solving $f(z)=0$. 

We let 
~$\displaystyle{K_f = \sum_{f(\zeta)=0} \log \frac{1}{\rho_\zeta}}$, 
and consider its average over all polynomials of a given degree
(including those with multiple roots).

\begin{thmC}
Let $\Lambar$ be the average value of $K_f/d$ over
$f\in\PDone$, where $\PDone$ is parameterized by the polydisk of
the roots endowed with Lebesgue measure.  
Then
$$\Lambar < 3d/2.$$
Consequently, the average of\/ $\Avgsteps{f}$ over $\PDone$ is $\BigO{d}$.
\end{thmC}

\begin{rem}\label{arith_complexity}
The cost of each step of the $\alpha$-step algorithm is dominated by the
calculation of $\alpha_f(z)$ (defined in Equation~(\ref{alpha})),
which can be done with  
$\BigO{d\log^2 d}$ arithmetic operations (see \cite{BorodinMunro}, for
example).
Consequently, Theorem~\ref{AverageCostBound} implies that for a specific
polynomial $f$, the expected arithmetic complexity to locate an approximate
zero via the $\alpha$-step algorithm is less than 
$\BigO{K_f \log^2d}$.
Combining this with
Theorem~\ref{LambdaLinearGrowth} gives an expected arithmetic complexity of  
$\BigO{d^2 \log^2d}$ to locate a root for a polynomial in $\PDone$.
\end{rem}

\begin{rem}\label{AllRootsRemarkIntro}
For $f\in\PDone$, by choosing $d$ appropriate starting values, an
approximate zero can be found for each root $\zeta_j$ in $\BigO{K_f}$ steps
of the $\alpha$-step algorithm.  This has an average arithmetic
complexity of $\BigO{d^3\log^2 d}$.  An explicit method for choosing initial
points is given in Section~\ref{AllRootsSec}.
\end{rem}

\smallskip
In addition to the theorems above, we wish to highlight several surprising
intermediate results which appear in Section~\ref{circlebounds}.
Specifically, let $|z_r|=r$ with $r>1$.  Then a bound on the rate of
change of $\Arg{f(z_r)}$ is given by our Angular Speed Lemma
(Lemma~\ref{argspeed}); applying this improves  
Proposition~2 of \cite{SS_Complexity2} regarding
the measure of ``good starting points'' from $1/6$ to $5/6$ (see 
Remark~\ref{speedremark}).  

Also worth noting are Corollary~\ref{Key1}, which shows that the average
value of $|f(z_r)|$ is $d\log{r}$, and Proposition~\ref{minabsoff}, which
states that  $|f(z_r)|$ is bounded below by a constant times $\rho_\zeta$.

\smallskip

\subsection*{Related Work}
 In \cite{Renegar}, Renegar gives an algorithm which approximates all
 $d$ roots of a polynomial with an arithmetic complexity of 
 $\BigO{d^3\log d + d^2\log d\log|\log\epsilon|}$ in the worst case.
 However, this algorithm includes a component requiring exact computation. 
 Pan's algorithm \cite{Pan} achieves the nearly optimal bound with a
 complexity of $\BigO{d^2\log d\log|\log\epsilon|}$, but
 implementation requires high precision computations (of the order exceeding
 the degree of the input  polynomial). 

 In practice, the software package MPSolve \cite{MPSolve} is widely
 used and empirical data indicates good global convergence
 properties;  the software uses the Aberth-Ehrlich method (see
 \cite{Aberth}, \cite{Ehrlich}) to locate the roots of the given
 polynomial.  There is not a lot of theoretical support, however: to our
 knowledge the global behavior of the Aberth-Ehrlich method is not
 understood. 
 
 In \cite{KS}, a worst-case complexity of $\BigO{d^2\log^2 d + d\log
   d|\log\epsilon|}$ yields an $\epsilon$-factorization for a polynomial $f$.
 This relies on a path-lifting algorthm which finds half the roots,
 then deflates the polynomial (that is, divides out by the
 approximations).  

 Recent work of Schleicher (\cite{SchFields}, \cite{BAS} and his co-authors
 have extended 
 the results of \cite{HSS} to obtain bounds for the complexity of finding
 $\epsilon$-roots. In \cite{HSS}, it is shown that there is a universal set
 of $1.1 d \log^2 d$
 points on a circle containing all the roots; if
 the roots are uniformly and independently distributed, \cite{BAS} shows
 that $\BigO{ d^2 \log^4 d}$ 
 iterations of Newton's method will locate all of the roots (an arithmetic
 complexity of $\BigO{d^3 \log^6 d + d^2 \log d\log|\log \epsilon|}$) with a
 high probability, comparable with the average arithmetic
 complexity of $\BigO{d^3 log^2 d + d^2\log d \log|\log \epsilon|}$ for the
 $\alpha$-step method in this paper (here the $\log|\log\epsilon|$
 term is added to account for the cost of refining an approximate zero to an
 $\epsilon$-root).

 One significant advantage of path-lifting methods over other methods
 is that of stability: as a consequence of estimates in \cite{K1}, as
 long as $f$ and its derivatives are computed with a relative error of
 $10^{-3}$, the algorithm will converge to an approximate zero in the
 same way. 


\subsection*{Organization}\label{intro_organization}
The paper is organized as follows.  In Section~\ref{prelim}, 
we set out notation and preliminary notions.  Section~\ref{alg}
describes the $\alpha$-step path-lifting algorithm explicitly.
In Section~\ref{branched}, we 
discuss the branched surface $\SS$ and the corresponding Voronoi partition.
This section discusses underlying topological and geometric
properties, and may be of interest independent to the question of
root-finding.   

Section~\ref{circlebounds} computes several estimates related to how
the polynomial $f$ behaves on the initial circle.
\Oldnote{mildly rephrased paragraph about \S\ref{optimize} and \S\ref{cost}}
In Section~\ref{optimize}, we bound the distance bewteen  $w_n$ and
$w_{n+1}$, and use this in \S\ref{cost} to estimate the number of
steps needed for the algorithm to  locate an approximate zero from a
given starting point~$z_0$, proving Theorem~\ref{PointwiseCostBound}.   

In Section~\ref{average}, we combine the topological and geometric
results of \S\ref{branched} with the more analytical results from
\S\ref{cost} to calculate an average upper bound over all starting
points for a given polynomial, proving Theorem~\ref{AverageCostBound}. 
In Section~\ref{CostVsDegree}, we discuss the relation between the
number of steps required and the degree of $f$ and proves
Theorem~\ref{LambdaLinearGrowth}.  This is followed by
Section~\ref{AllRootsSec} where we describe how to use this method to
locate all roots of a polynomial $f$.  We conclude in Section~\ref{Remarks}
with some remarks and comments regarding extensions of these results.

\subsection*{Acknowledgement}
The authors would like to thank Araceli Bonifant, Michael Shub, and the
anonymous referees for their input and suggestions which have improved this
paper.

\section {Preliminaries}
\label{prelim}
We will use the following general notions and notations throughout.
\smallskip

An open disk of radius $r>0$ centered around $z\in \C$ is denoted by 
$D_r(z)$.  

Let $S_r(z)$ denote the circle of radius $r$ and center $z$; if the 
circle is centered at the origin, we will denote it by $S_r$.

The function $\Arg$ denotes the argument of a complex number 
(in the interval $(-\pi,\pi]$ unless otherwise noted).

The \dfn{ray} $\ray{w}\subset \C$ of a point $w\in \C\smallsetminus\{0\}$ is   
$$\ray{w}=(0,\infty)\cdot w = \Set{ z\in\C \st \Arg{w} = \Arg{z}},$$
and the \dfn{slit} of this point is the part of the ray extending
outward from $w$, 
that is 
$$\slit{w}=[1,\infty)\cdot w = \Set{ z\in \ray{w} \st |z| \ge |w|}$$ 

For a polynomial $f:\C\to \C$, denote the critical points of $f$ by
$$
\CC_f=\Set{z \st f'(z)=0}.
$$ 
For a regular point $z_0$, we shall use $f^{-1}_{z_0}$ to denote a holomorphic
branch of 
the inverse of $f$ for which $f^{-1}(f(z_0)) = z_0$.

\begin{figure}[htpb]
\begin{center}
\begin{overpic}[width=.5\hsize,tics=5]{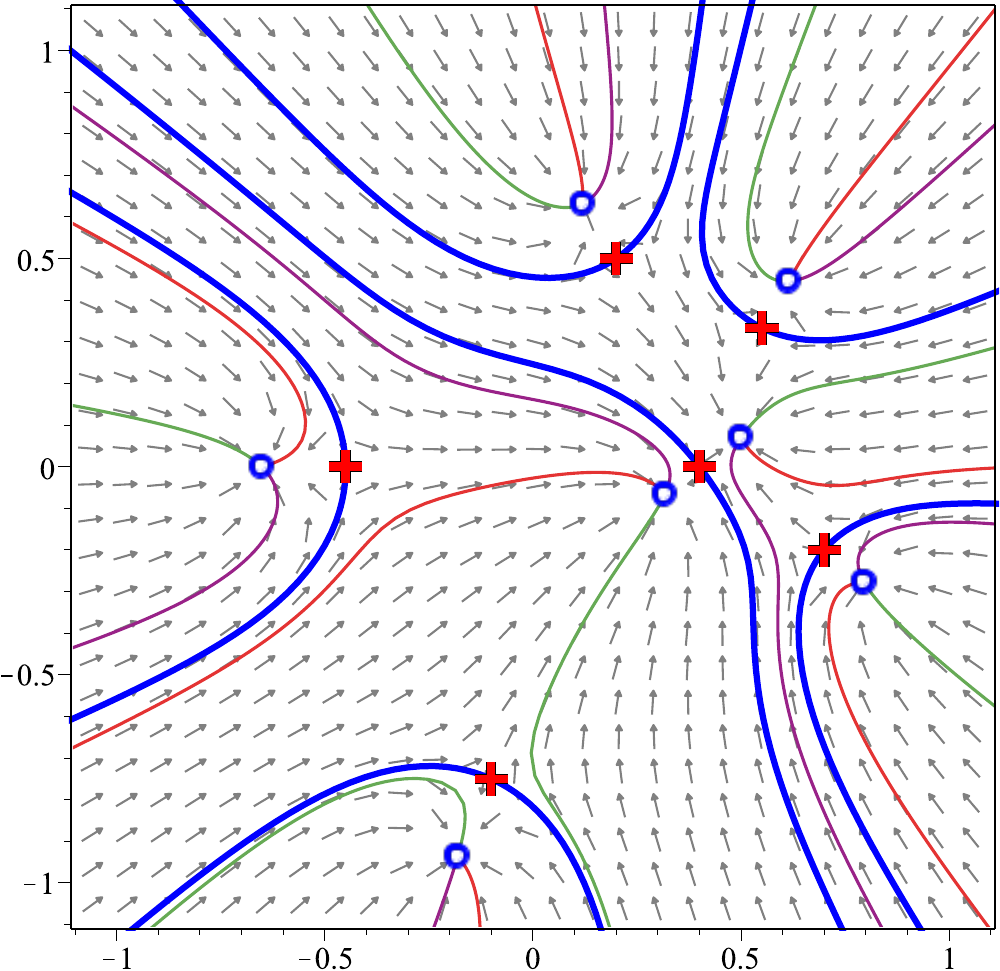}
\end{overpic}
\OldnoteC{Referee found orig fig~\ref{newtflowpic} ``hard to digest''; This
        is redone. Better? (I removed labels on $c_j$ and $\zeta_k$ because it
      was too cluttered.)}   
\mycaption{\label{newtflowpic} The direction field for the Newton flow
  corresponding to a degree~7 polynomial is shown.  For each
  root $\zeta_i$ (indicated by a circle \CIRC),
 its basin is bounded by the stable manifolds (thick curves \LINE) 
 of one or more critical points $c_j$ (indicated by a cross \CROSS). 
 Also shown are solution curves $\varphi(t)$ for which $\Arg{f(\varphi)}$ is
 $0$, $2\pi/3$, or $-2\pi/3$ (thin curves).
Compare Figure~\ref{RSpicture}.
}
\end{center}
\end{figure}

\medskip
\label{NewtFlow}
We now discuss the Newton flow, and some notation related to it.
Consider the following vector field on $\C$,
$$
X(z)=-\frac{f(z)}{f'(z)}.
$$
The corresponding flow is called the \dfn{Newton flow}. This vector
field blows up near the critical points of $f$. By rescaling the
length of the vector $X(z)$ by $2|f'(z)|^2$,  the critical points of
$f$ become well-defined singular points of the rescaled vector
field. This rescaled vector field is the gradient vector field
                  $\dot{z} = -\nabla |f(z)|^2$; the
solution curves of the former coincide with the
latter, and we will use the two interchangably.  The equilibria of
the Newton flow are exactly the roots and critical points of $f$.
Each root $\zeta$ is a sink; we shall denote its basin of attraction
by $\Bas(\zeta)$. Critical points are saddles for the flow.
Furthermore, we can extend the flow to infinity, which is the only
source. Each boundary component of $\Bas(\zeta)$ contains critical
points $c\in \CC_f$: generically, each critical point $c$ has an
unstable orbit leaving from $c$ and converging to $\zeta$, as well as 
stable orbits from infinity to $c$, which are separatrices for the flow.
Generically, there is a unique critical point in each boundary
component; in the degenerate cases, there could be saddle connections
resulting in multiple critical points on one boundary component.
A general discussion regarding the Newton flow can be found in \cite{STW}
and~\cite{JJT}, as well as \cite{KozenStefansson}.   See
Figure~\ref{newtflowpic}.  

It is important to note that if $\varphi(t)$ is a solution curve for the
Newton flow, $f(\varphi(t))$ lies along a ray. 
 To see this, observe that
$$\frac{d}{dt}f(\varphi(t)) = 
  f'(\varphi(t))\cdot \left( -\frac{f(\varphi(t))}{f'(\varphi(t))}\right) =
  -f(\varphi(t)) ,$$
and hence $f(\varphi(t)) = e^{-t}f(z_0)$ for some $z_0=\varphi(t_0)$,
provided $f'(\varphi(t))$ is never zero. (If $\varphi(t)$ contains a critical
point of $f$, the result follows by continuity.)

Since $f$ has distinct roots, $f'(\zeta)\ne 0$ for each root $\zeta$, and so
$f$ is a local diffeomorphism in a neighborhood of $\zeta$. Thus, for every
angle $\theta$ there will be a solution $\varphi_\theta(t)$ in this
neighborhood with $\Arg(f(\varphi_{\theta}(t)))=\theta$.  Noting that the ray
$f(\varphi_\theta(t))$ extends to infinity unless $\varphi_\theta(t)$
encounters a critical point $c$, we obtain the following lemma.

\begin{lem}\label{EachBasinCoversPlane}
For each root $\zeta$, $f$ is a biholomorphism
$$f\colon\Bas(\zeta) \rightarrow \C \smallsetminus \bigcup \slit{f(c)} ,$$
where the union is taken over the critical points $c$ which lie on the
boundary of $\Bas(\zeta)$.
\end{lem}

\smallskip
\begin{rem}\label{NewtonAsEuler}
Observe that iteration of Newton's method beginning at a point $z_0$
corresponds to construction of an approximate solution to the Newton flow with
intial condition $\varphi(0)=z_0$ using Euler's method with stepsize $h=1$.
When the path $\ThePath$ is a ray in the target space, a 
path-lifting method corresponds to constructing approximate solutions of the
Newton flow via a method that self-corrects to always follow a solution curve
that containing the initial condition.
\end{rem}

\bigskip
Throughout the paper, we will consider polynomials $f \in \PDone$,
that is, $f:\C\to \C$ given by
$$f(z)=\prod_{j=1}^d (z-\zeta_j) \qquad\mbox{with}\ |\zeta_j|\le 1,$$
with distinct roots $\zeta_j$.
 The set of roots of $f$ will be denoted by 
 $$\Roots_f=\Set{\zeta_j \st j=1,\dots,d}.$$

The restriction to $\PDone$ is not severe; provided its roots are simple, an
affine change of coordinates depending only on the coefficients 
\Oldnote{added ``simple roots'' here.}
will transform any polynomial into one in $\PDone$ (see \cite{GeomPolys},
for example).  
%
%
The space $\PDone$ is somewhat different from that considered in other
works (such as 
$\mathbf{P}_1$ of \cite{Sm81}, ${\mathcal P}{\kern-0.2em}_d(1)$ of
\cite{KS}, etc.), where the space of 
polynomials is 
represented
as $\Set{f(z)=\sum a_jz^j \st  |a_j| \le 1}$.  
In this case, all the roots lie in the disk of radius $2$, and
our results are readily adapted to any set of polynomials where the
roots lie in any disk of a known radius.

\smallskip
We shall use the following standard result several times.

\begin{lem}[{\bfseries Koebe Distortion Theorem}]\label{Koebe_Lemma}
Let $g: D_r(0) \to \C$ be univalent with 
$g(0)=0$ and $g'(0)=1$. For $z\in D_r(0)$ with $s=|z|/r$, we have
\begin{equation}\label{KLdiff}
\frac{1-s}{(1+s)^3}\le |g'(z)|\le \frac{1+s}{(1-s)^3}
\end{equation}
and  
\begin{equation}\label{KLzbound} 
 \frac{|z|}{(1+s)^2}\le |g(z)|\le \frac{|z|}{(1-s)^2}
\end{equation}
Consequently,
\begin{equation}\label{Koebe14}
D_{r/4}(0)\subset g(D_r(0)).
\end{equation}
\end{lem}

\begin{rem} The statement in \EqRef{Koebe14} is known as the
  Koebe~$\frac14$-Theorem. The proof can be found in \cite{Ko},
  \cite{P}, or \cite{Duren}, among others. See also Corollary~2.6 of \cite{K2}.
\end{rem}
\bigbreak

\section{The Path-Lifting Algorithm}\label{alg}
In this section,  we present the  path-lifting algorithm that we use to 
find an approximate zero of $f \in \PDone$.  First, we discuss approximate
zeros and the Kim-Smale $\alpha$~function.
\Oldnote{changed $z_0$ to $z^*$ in Def.~\ref{approxzero}; added statement about
quadratic convergence at end}

\begin{defn}\label{approxzero} Let $z_n\in \C$ be the $n\Th$
iterate under Newton's method of the point $z^*\in \C$, that is, 
$$
z_{n+1}=z_n-\frac{f(z_n)}{f'(z_n)}, \quad z_0=z^* .
$$ 
The point $z_*$ is called an \dfn{approximate zero} of $f$ if 
$$
|z_{n+1}-z_n| \le \left(\frac{1}{2}\right)^{2^n-1} |z_1-z^*|
\qquad\text{for all $n>0$.}  
$$ 
\Oldrnote{Def~\ref{approxzero}: Referee asks  ``is this ever possible exactly,
  not just asymptotically?'' 
  I have no idea what he means.  Does he mean ``is the
  inequality ever sharp for all $n$?'' or ``do approximate zeros exist?''}
Newton's method converges quadratically to a root when started from an
approximate zero (see \cite{Sm85} for example). 
\end{defn}

Approximate zeros are an intrinsic, dynamical feature of a polynomial.
They form disjoint connected neighborhoods of the roots $\zeta_i$ on which the
Newton map $N_f(z)=z-f(z)/f'(z)$ converges quadratically to the root,
which is a super-attracting fixed point for the rational map $N_f$.

A sufficient condition for a point to be an approximate zero is developed in 
\cite{K1} and \cite{SmaleDataAtOnePoint}. We will use the criterion formulated by Smale in
\cite{SmaleDataAtOnePoint} to locate approximate zeros. It uses
$\alpha:\C\smallsetminus \CC_f\to\R$ defined by 
\begin{equation}\label{alpha}
\alpha(z)=\max_{j>1}\biggl|\frac{f(z)}{f'(z)}\biggr| 
\biggl|\frac{f^{(j)}(z)}{j!f'(z)}\biggr|^{\frac{1}{j-1}}.
\end{equation}
It is sometimes useful to use the related function $\gamma(z)$
instead, where 
\begin{equation}\label{gamma}
\gamma(z)=\max_{j>1} \left|\frac{f^{(j)}(z)}{j!f'(z)}\right|^{\frac{1}{j-1}}.
\end{equation}
While we will primarily use $\alpha(z)$, we make use of $\gamma(z)$ in
Corollary~\ref{AlphaOutsideDisk}, Section~\ref{optimize} and
Section~\ref{CostVsDegree}.

\medskip
\begin{thm}\label{alphaapproxzero} (\cite{K1},\cite{SmaleDataAtOnePoint})
There is a number $\alpha_0$ such that if 
$\alpha(z) < \alpha_0$, the point $z$ is an approximate zero.
\end{thm}

\begin{rem}
It has been shown that $\alpha_0 \ge 3-\sqrt{8} \approx 0.17157$ 
(see \cite{WH} or \cite{WZ}, for example).
\end{rem}

\begin{rem}
The number $\alpha_0$ is given in \cite{SmaleDataAtOnePoint} and in many places throughout
the literature as $\alpha_0 \approx 0.130707$.  However, this specific value
is very likely the result of a typographic error in the fifth decimal place.
Smale's bound for $\alpha_0$ is stated as a solution to 
$(2r^2 - 4r +1)^2 - 2r = 0$ \cite[Section~4] {SmaleDataAtOnePoint}; 
the relevant root of this equation is $0.130716944\ldots$.   
\end{rem}

\bigbreak
We shall analyze the following algorithm to find an approximate zero for 
$f \in \PDone$. 

\bigskip

\framebox{
\begin{minipage}{.8\textwidth}
\centerline{\bf The $\alpha$-Step Path-Lifting Algorithm}\label{plm}
\medskip
Input a polynomial $f\in\PDone$. 
\OldnoteC{added ``Input $f$'' to statement of algorithm, and output $z_n$}
\begin{description}
\item[{\bf Step 0}] Choose $z_0\in \C$ with $|z_0|=1+\frac1d$. 
Let
$$
w_0=f(z_0)  \quad\hbox{and}\quad  w=\frac{w_0}{|w_0|}.
$$
\item[{\bf Step 1}] Stop if $\alpha(z_n)\le \alphaZero$;\ Output
  $z_n$,  an approximate zero for $f$. \\

\item[{\bf Step 2}] Let
$$
w_{n+1}=w_n-\frac{1}{15}\cdot \frac{|f(z_n)|}{\alpha(z_n)}\cdot w
$$
and 
$$
z_{n+1}=z_n-\frac{f(z_n)-w_{n+1}}{f'(z_n)}.
$$
Continue with Step 1.
\end{description}
\end{minipage}
}

\medskip\noindent
Sometimes we shall refer to the points $w_n$ generated by the algorithm
above as \dfn{guide points} or \dfn{target points}.
\bigbreak

If $z_0\in \Bas(\zeta)$ then the $\alpha$-step algorithm will terminate with an
approximate zero for $\zeta$.
\Oldnote{expanded on algorithm terminating}
This follows from the fact that $\Arg{w_n} = \Arg{f(z_0)}$ for all $n$, and,
by the estimates in Section~\ref{optimize}, 
$w_{n+1}$ is close enough
to $w_n$ to ensure that $f^{-1}$ is univalent on a region containing $w_n$,
$w_{n+1}$, $f(z_n)$, and~$f(z_{n+1})$.
Since $z_0 \in \Bas(\zeta)$, the entire ray $\ray{w_0}$ lifts to a curve
lying in $\Bas(\zeta)$ since~$\ray{w_0}$ does not contain a critical
value $f(c)$ with $c$ in the closure of $\Bas(\zeta)$.

\begin{figure}[tbhp]
\begin{center}
\begin{overpic}[width=.9\hsize,tics=5]{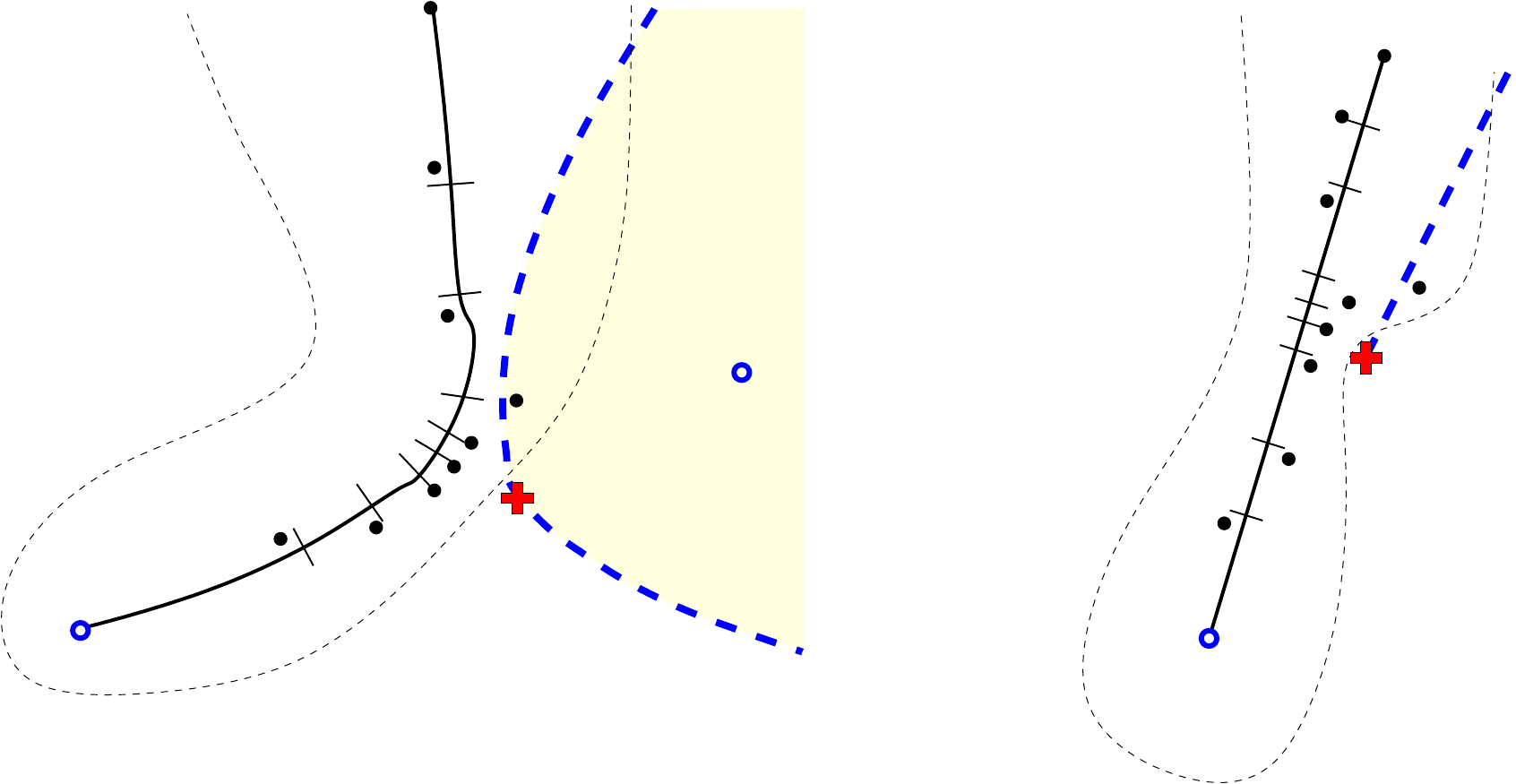}
\put(4,7)  {$\zeta_1$}
\put(48,24){$\zeta_2$}
\put(80,7) {$0$}
\put(30,5) {\small$\Bas(\zeta_1)$}
\put(44,15) {\small$\Bas(\zeta_2)$}
\put(10,10){\small $f^{-1}_{\zeta_1}(\ThePath)$}
\put(82,12){$\ThePath$}
\put(33.5,15.5){$c$}
\put(91,27){\small $f(c)$}
\put(83,49){$U$}
\put(14,49){$f^{-1}_{\zeta_1}(U)$}
\put(29,51){\tiny $z_0$}
\put(26,40){\tiny $z_1$}
\put(27,30.5){\tiny $z_2$}
\put(35,25){\tiny $z_3$}
\put(18,17.2){\tiny $z_N$}
\put(87,49.5){\tiny $f(z_0)=w_0$}
\put(84,44)  {\tiny $f(z_1)$}
\put(91,43)  {\tiny $w_1$}
\put(82.5,38){\tiny $f(z_2)$}
\put(90,39)  {\tiny $w_2$}
\put(94.5,32){\tiny $f(z_3)$}
\put(84,34)  {\tiny $w_3$}
\put(83,32)  {\tiny $w_4$}
\put(84,17)  {\tiny $w_N$}
\put(75.5,16.5){\tiny $f(z_N)$}
\put(50,35){$\displaystyle \xymatrix{\relax& \ar@/^{.8pc}/[rr]_{f} &\relax&\relax}$}
\end{overpic}
\OldnoteC{Figure \ref{plmdef} redone to agree with conventions in other figures}
\mycaption{An illustration of the $\alpha$-step method beginning at $z_0$,
  with $\C_{source}$ on the left and $\C_{target}$ on the right.  The guide
  points $w_j$ (and their preimages) are shown along $\ThePath$ and
  $f_{\zeta_1}^{-1}(\ThePath)$ as the intersection of perpendicular
  segments. The points $z_0$ and their images $f(z_0)$ are indicated by
  solid dots, two roots $\zeta_1$ and $\zeta_2$ (and their image $0$) are
  denoted by circles~(\CIRC), and a nearby critical point $c$ and its image
  $f(v)$ are marked by a cross~(\CROSS). $\Bas(\zeta_2)$ is shaded.
  In this illustration, $z_0 \in \Bas(\zeta_1)$ but $z_3\in\Bas(\zeta_2)$.
  However, as noted in Remark~\ref{znotinbasin}, there is a neighborhood $U$
  of the ray on which there is a branch of the 
  inverse which contains all the $z_n$.  $U$ is shown bounded by a
  dashed line.}
 \label{plmdef}
\end{center}
\end{figure}

\begin{rem}\label{znotinbasin}  There may be some values of $n$
for which $z_n\notin \Bas(\zeta)$; even if this occurs, 
there is a neighborhood $U\subset \C$  of the ray $\ray{w_0}$ which
contains $f(z_j)$ for all $j$ and on which there exists a univalent
inverse branch of $f^{-1}$ mapping $w_0$ to $z_0$.   
\Oldnote{Referee asked where is Rem.\ref{znotinbasin} shown. added pointer to
  \S\ref{ConstructingU}} 
As noted in the previous paragraph, $w_{n+1}$ and~$f(z_n)$ both lie in
a neighborhood of $\ray{w_0}$ on which $f^{-1}$ is univalent, even if
$z_n$ is outside $\Bas(\zeta)$. In this case, $\Bas(\zeta)$ can be
enlarged to a neighborhood $U$ of $\ThePath$ which contains all the $z_j$.
See Figure~\ref{plmdef}.  A more detailed description and explicit
construction of $U$ can be found on page~\pageref{ConstructingU} of 
Section~\ref{ConstructingU}. 
Denote this inverse branch by $f^{-1}_{z_0}:U\to\C$. 
\end{rem}


\begin{defn}\label{rhoDef}
For every zero $\zeta\in \Roots_f$, let
\[
\rho_\zeta= \min_{c\in \CC_f(\zeta)} |f(c)|
\qquad\text{where}\qquad
\CC_f(\zeta)=\CC_f\cap \overline{\Bas(\zeta)}.
\]
\end{defn}

\begin{rem}\label{RhoIsDistToV} 
Note that $\rho_\zeta$ is the radius of convergence of $f^{-1}_\zeta$
at $0$, and is the distance in the surface $\SS$ between $\fhat(\zeta)$
and the nearest branch point of $\SS$. 
\Oldnote{added justification of rem~\ref{RhoIsDistToV}, since referee
  doubted it.}
This follows from the fact that 
$\fhat\colon \Bas(\zeta) \rightarrow \SS \smallsetminus \VV_f$
is a biholomorphism 
and $\pi$ is an isometry (see Lemma~\ref{Ball}) 
from the disk~$D_{\rho_\zeta}$ about $\fhat(\zeta)$ into $\C_{target}$. 
Hence, $f^{-1}_\zeta \colon D_{\rho_\zeta}(0) \rightarrow \C_{source}$ is a
univalent analytic function.
\end{rem}

\begin{defn}\label{Kfdef}
For any polynomial $f$, we define
$\displaystyle{
K_f= \sum_{\zeta \in \Roots_f} \log \frac{1}{\rho_\zeta}
}$.
\end{defn}

\begin{rem}\label{KfInfinite} Notice that $K_f<\infty$ if and only if 
the set of roots $\Roots_f$ and critical points $\CC_f$ are disjoint.
 This holds generically for polynomials $f$, and
  $K_f=\infty$ exactly when $f$ has a multiple zero.  Root-finding problems for
  which there is a multiple zero are typically called \dfn{ill-conditioned} or 
  \dfn{ill-posed}.
\end{rem}
\Oldnote{removed former remark~3.11 that called $K_f$ a condition number,
  since it isn't clear what metric is appropriate here, and it made an 
  earlier referee mad (inisting the THE condition number for
  root-finding was $\mu_f$). Not worth the trouble.}

\begin{rem} One can introduce a 
measure of difficulty
$K_{f,\zeta} = \log 1/\rho_\zeta$ for a specific given root $\zeta\in \Roots_f$. 
Then Theorem~\ref{PointwiseCostBound} describes the cost of reaching
an approximate zero for $\zeta$ in terms of $K_{f,\zeta}$, 
Theorem~\ref{AverageCostBound} gives the cost of finding any
approximate zero in terms of the average value of $K_{f,\zeta}$, 
and Theorem~\ref{LambdaLinearGrowth} averages $K_{f,\zeta}$ over all
polynomials $f$ of a given degree.
\end{rem}
\Oldrnote{removed remark about $\alpha$ when $f$ is an iteration. Also
removed the end of this section, since it was a rehash of ``Main Results''
in introduction.}

\section{The Voronoi Partition in the Branched Cover}\label{branched}

Given a polynomial $f:\C\rightarrow \C$ of degree $d$, recall 
from Section~\ref{prelim} \Oldnote{added ref. for ``recall''}
that we denote its
critical points by $\CC_f=\Set{z \st f'(z)=0}$.  For any such $f$, we
can express it as a composition $f=\pi\circ\fhat$,
$$  \xymatrix{
  \C_{source}    \ar[r]^{\fhat}  
        \ar[dr]_{f}     & \SS  \ar[d]^{\pi}\\
  \relax                & \C_{target} },
$$
\Oldnote{added ``source'' and ``target'' to diagram, changed diffeo to biholo,
  added ``see also''}
where $\fhat$ is a biholomorphism except on $\CC_f$ 
(on which it is merely a bijection), and $\pi$ is a $d$-fold branched cover,
ramified at points of $\VV_f=\fhat(\CC_f)$.   We deonte the metric on
$\SS$ by  $\dist{\cdot,\cdot}$; this metric is such that away from points in
\Oldnote{revised sentence about metric on $\SS$ due to ref's comment.}
$\VV_f$, $\pi$ is a local isometry into $\C_{target}$ (with the standard
metric).  See also Figure~\ref{RSpicture} and the corresponding discussion
in Section~\ref{intro_background}.  

\medskip
The \dfn{multiplicity} of a critical point $c\in \CC_f$ is
$$
m_c=\min \Set{k \st f^{(k+1)}(c)\ne 0}. 
$$
Notice that
$$
\sum_{c\in \CC_f} m_c =d-1.
$$
The points in $\VV_f$ are called \dfn{critical values} in $\SS$, and
we define the \dfn{multiplicity} $m_v$ of $v=\fhat(c)\in \VV_f$  to be the 
multiplicity of $c$;
\Oldnote{added phrase about degree of ramification}
this is also the local degree of the projection $\pi$ in a neighborhood of
$v$.   

\medskip
Note that for each root $\zeta \in \Roots_f$, 
$$\pi: \fhat\left(\Bas(\zeta)\right) \to 
       \C\smallsetminus 
       \kern-.3em\bigcup_{y\in V_f(\zeta)}\kern-.3em \slit{y}$$
is an isometry (where $V_f(\zeta) = f(\CC_f(\zeta))$, and $\slit{y}$ is
the ray outward from $y$ as defined in Section~\ref{prelim}).
\medskip

The \dfn{Voronoi domain} of a point $v\in \VV_f$ is 
$$
\Vor{v}=\Set{u\in \SS \st \dist{u,v} \le \dist{u,w}, \forall w\in
\VV_f};
$$ 
\Oldrnote{rephrased phrase after def of $\Vor{v}$; included def of ``relevant''
  per referee.  OK?}
this is exactly the set of points $u\in\SS$ such that the critical value
$\pi(v)$ lies on the boundary of the disk about $\pi(u)$ on which the
inverse $f^{-1}_x$ will be analytic ($x$ satifies $\fhat(x)=u$).  See
also Remark~\ref{RhoIsDistToV}.  We will refer to such critical values $\pi(v)$
as \dfn{relevant} to the construction of $f_x^{-1}$.

%
Recall 
from Section~\ref{prelim} \Oldnote{added ref. for ``recall''; replaced Euclidean}
that $D_r(u) = \Set{ y \st \dist{u,y} < r}$ denotes the open disk of
radius $r$ about $u$.  For $u\in\SS$, such disks will 
be isometric to their projections (i.e., be ``Euclidean disks'')
exactly when they avoid the branch points of $\SS$.

\begin{lem}\label{Ball} A point $u\in\SS$ is in $\Vor{v}$  if and only if  
$\pi: {D_{\dist{u,v}}(u)}\to  {D_{|u-v|}(\pi(u))}$
is an isometry. 
\Oldrnote{ref asks if issue because $\Vor{v}$ is closed.  Not a problem, I think}
In particular, if $u\in \Vor{v}$ then
$$
D_{\dist{u,v}}(u)\cap \VV_f=\emptyset.
$$
\end{lem}

\begin{proof}
 If $u\in \Vor{v}$ then $D_{\dist{u,v}}(u)\cap \VV_f=\emptyset$. 
Thus, $\pi$ is a local isometry on all of $D_{\dist{u,v}}(u)$, and in
particular,  
$\pi$ is a global isometry on this disk.  Conversely, 
If $\pi$ is an isometry on all of $D_{\dist{u,v}}(u)$, there can be no
critical values in the disk, and so $u\in \Vor{v}$. 
\end{proof}

Let $u_1, u_2\in \SS$. If the line segment $[\pi(u_1),\pi(u_2)]\subset\C$ 
has a lift in $\SS$ which connects $u_1$ with $u_2$, we denote this lifted
line segment by $\LiftSeg{u_1,u_2}$. Observe that many pairs $u_1,u_2$ do not
have such a connecting line segment. In this case we write
$\LiftSeg{u_1,u_2}=\emptyset$.   When $\LiftSeg{u_1,u_2}$ is nonempty, we say
that \dfn{$u_1$ is visible from $u_2$} in $\SS$.
Also observe, if $v \in \VV_f$ then
$$
\LiftSeg{u,v}\ne \emptyset \quad\text{for all}\ u\in \Vor{v}.
$$

\medskip

\begin{figure}[htpb]
\begin{center}
\begin{overpic}[width=.5\hsize,tics=5]{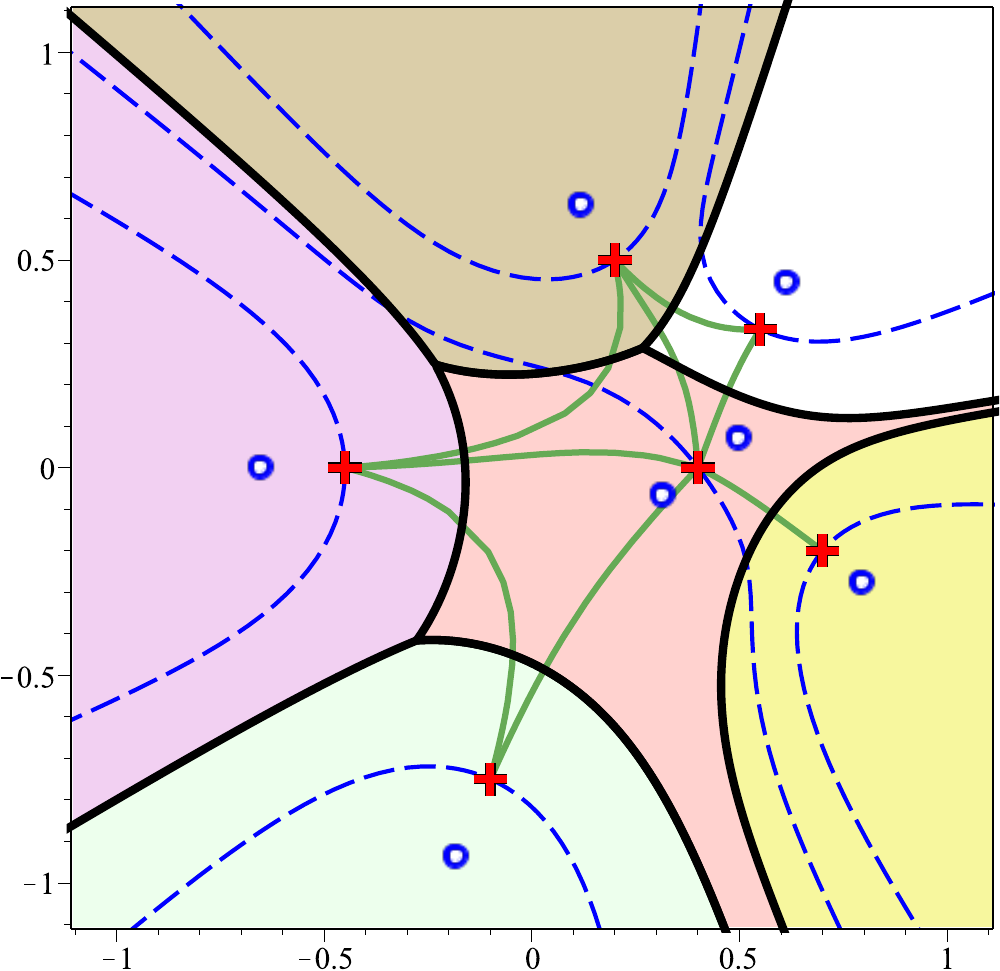}
\put(35,53) {\small $c_4$}
\put(51,19) {\small $c_5$}
\put(64,70) {\small $c_3$}
\put(77,60.5) {\small $c_2$}
\put(81,46) {\small $c_6$}
\put(72,49) {\small $c_1$}
\put(29,7)  {\small $\Vor{c_5}$}
\put(57,30) {\small $\Vor{c_1}$}
\put(85,20) {\small $\Vor{c_6}$}
\put(85,85) {\small $\Vor{c_2}$}
\put(40,85) {\small $\Vor{c_3}$}
\put(10,40) {\small $\Vor{c_4}$}
\end{overpic}
\OldnoteC{revised Fig~\ref{voronoipic} and caption}
\mycaption{\label{voronoipic} The Voronoi regions of Fig.~\ref{vorsheets}
  are shown in the source space $\C_{source}$.
  The roots of $f$ are indicated  by circles (\CIRC), the critical points by
  crosses (\CROSS) and labeled as $c_j$.
  The dashed lines are the boundaries of $\Bas(\zeta_j)$ for each root;
  each such boundary contains a unique critical point $c_k$; observe that
  each Voronoi domain enters the basin of at least two roots.
  For each critical point $c_j\in\CC_f$, $\fhat^{-1}(\Vor{v_j})$ is
  shown bounded by the heavy solid lines, shaded as in Fig.~\ref{vorsheets},
  and labeled as $\Vor{c_j}$.  
  The visibility graph $\fhat^{-1}(\VisG)$ is also shown, indicated by
  solid curves connecting pairs of critical points $c_j$ and $c_k$.
}
\end{center}
\end{figure}

We can form the \dfn{visibility graph for $\SS$} as follows.
\Oldnote{changed ``homeomorphism'' to ``bijection'' (not a homeo on $\CC_f$) regarding $\VisG$}
The vertices of the graph are the critical values $\VV_f$, and there
is an edge from $v$ to $w$ if and only if $\LiftSeg{v,w}$ is
non-empty. 
We can identify the visibility graph with the subset of $\SS$ given by 
$$\VisG = \bigcup_{v,w \in \VV_f} \LiftSeg{v,w}.$$
Since $\fhat$ is a bijection between $\C_{source}$ and $\SS$,
$\fhat^{-1}(\VisG)$ is well-defined, so we can also 
view $\VisG$ as a graph immersed in $\C_{source}$, with the critical points
of $f$ as vertices.

\begin{problem}\label{WhichGraphsAreVisGraph}
Characterize the graphs which occur as a visibility graph $\VisG$ for a
polynomial. 
\end{problem}
\Oldnote{reworded question\ref{WhichGraphsAreVisGraph}; 
   Removed 5 paragraphs about $L_{u,v}$, $\NN_v$ following
   it, as well as a figure.}

Recall 
from Section~\ref{prelim} 
that the  ray $\ray{y}\subset \C$ of a point $y\in
\C\smallsetminus\{0\}$ is  the set of points which have the same
argument as $y$. 

If $\widehat{0}\in \SS$ projects onto $0$ and 
$\LiftSeg{\widehat{0},u}\ne \emptyset$, the 
geodesic starting at $\widehat{0}$ and containing $\LiftSeg{\widehat{0},u}$
is the ray through $u\in \SS$, which we denote by ${\hatray{u}}$. 
Observe that if $\hatray{u}\cap \VV_f=\emptyset$ then
$
\pi: {\hatray{u}}\to \ray{\pi(u)}
$
is a surjective isometry.

Let $y=\pi(u)$.  If $\ray{y}\cap f(\CC_f)=\emptyset$, then 
$$
\pi^{-1}(\ray{y})={\hatray{y_1}}\cup {\hatray{y_2}}\cup
\dots \cup {\hatray{y_d}},
$$
where the points $y_i\in \SS$ are the $d$ different preimages of $y$.

\begin{prop}\label{4pie} Given $v\in \VV_f$ and 
 $y\in \C\smallsetminus  f(\CC_f)$. 
 Then
$$
\card\Set{i \st {\hatray{y_i}}\cap \Vor{v}\ne \emptyset }\le m_v+1.
$$
Furthermore, each $\hatray{y_i}\cap \Vor{v}$ is a connected set.
\end{prop}

\begin{proof} Suppose ${\hatray{y_1}}, {\hatray{y_2}},\dots, 
{\hatray{y_k}}$ intersect $\Vor{v}$, with $v=\fhat(c)$, $c\in\CC_f$.
 Pick a point $u_i$ in each
of these intersections, that is, 
$$
u_i\in {\hatray{y_i}}\cap \Vor{v}.
$$

\begin{figure}[htbp]
\begin{center}
\begin{overpic}[width=0.6\hsize,tics=5]{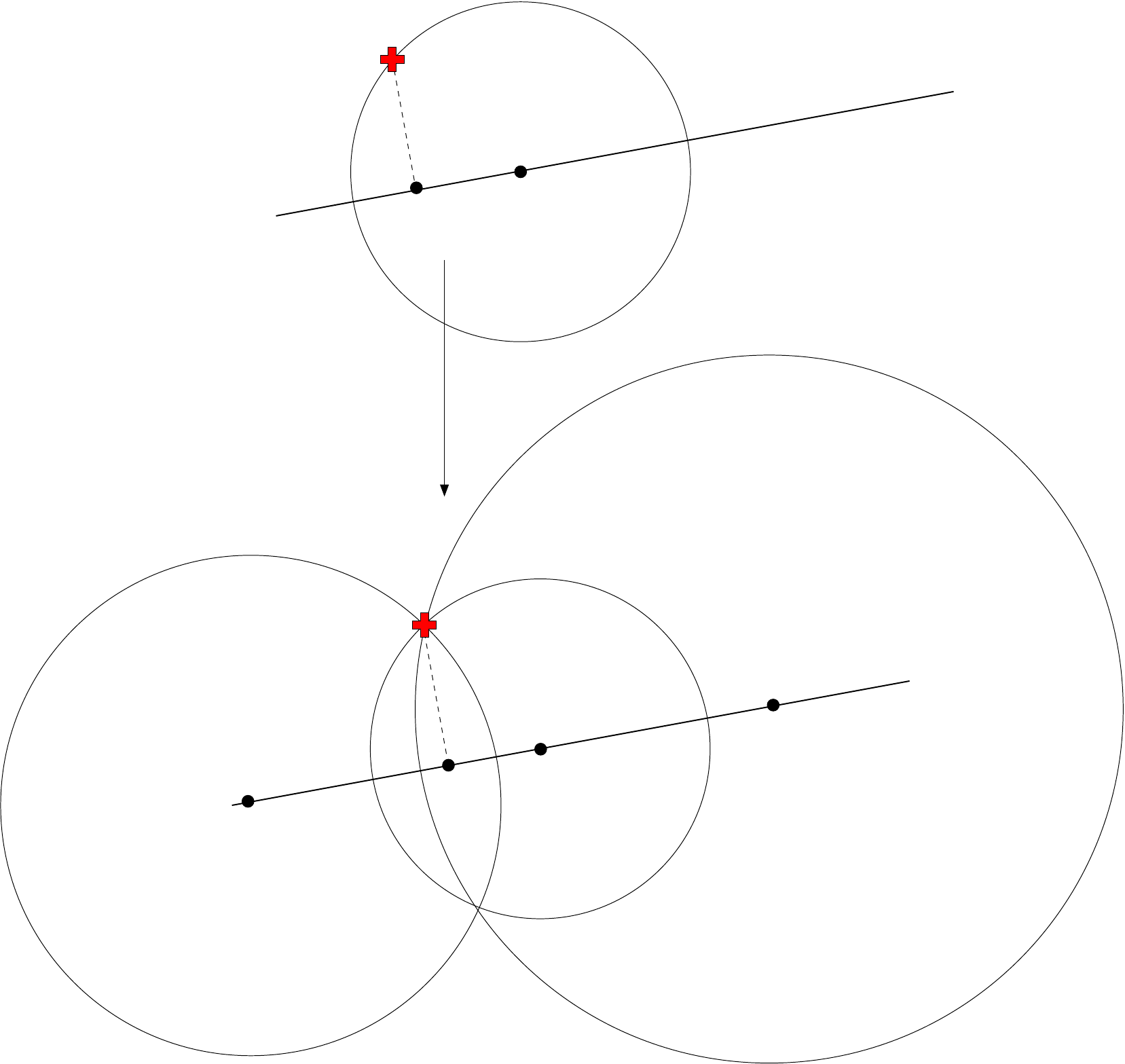}
\put (23,91) {\small $\fhat(c)=v$}
\put (29,38) {\small $\pi(v)$}
\put (36,60) {\small $\pi$}
\put ( 6,35) {\small $\pi(D_1)$}
\put (50,69) {\small $D_i$}
\put (50,18) {\small $\pi(D_i)$}
\put (80,50) {\small $\pi(D_k)$}
\put (45,81) {\small $u_i$}
\put (36,75) {\small $p_i$}
\put (39,23) {\small $p$}
\put (46,24) {\small $\pi(u_i)$}
\put (20,20) {\small $\pi(u_1)$}
\put (67,28) {\small $\pi(u_k)$}
\put (81,33) {\small $\ray{y}$}
\put (85,85) {\small $\hatray{y}$}
\end{overpic}
\OldnoteC{redid fig~\ref{prop} to use \CROSS\ for critical value}
\mycaption{As proven in Proposition~\ref{4pie}, the
  projection $\pi$ is at most $(m_v+1)$-to-one on $\Vor{v}$.} \label{prop} 
\end{center}
\end{figure}
Let $D_i=D_{\dist{v,u_i}}(u_i)$. According to Lemma~\ref{Ball}, we know that
$\pi: D_i\to \pi(D_i)$
is an isometry. Let $p_i\in {\hatray{y_i}}$ be the perpendicular
projection of $v$ onto ${\hatray{y_i}}$ and let $p$ be the
projection of $f(c) = \pi(v)$ onto~$\ray{y}$. See Figure~\ref{prop}
Then for all $i\le k$, $\pi(p_i)=p$, 
$$
\emptyset\ne \LiftSeg{v,p_i}\subset D_i
\qquad\text{and}\qquad
\emptyset\ne [\pi(v),p]\subset \bigcap_{i\le k} \pi(D_i).
$$
Since for each $i$ between $1$ and $k$, $\pi$ is a surjective isometry from 
$\LiftSeg{v,p_i}$ to $[\pi(v),p]$, $k$ can be no larger than the the degree 
of $\pi$ on a neighborhood of $v$.  That is,
  $$k \le 1 + m_v .$$

The connectedness of $\hatray{y_i}\cap\Vor{v}$ follows from the triangle
inequality.
\end{proof}

\begin{cor}\label{twotoone} Each projection $\pi:\Vor{v}\to \C$ is at
  most $(m_v+1)$-to-one.  
\end{cor}


\bigskip
Let $z\in\C$.  We'll say that a critical point $c\in\CC_f$
\dfn{influences the orbit of $z$} if the segment
$\LiftSeg{\widehat{0},\fhat(z)}$ passes through $\Vor{\fhat(c)}$.
\medskip 

We are interested in the critical points which influence the starting
points for our algorithm, and, conversely, the starting points which
are influenced by a given critical point.  

\begin{defn}\label{InflSets}
For starting points $z$ on the
circle of radius $r$, we define the following sets:
\begin{equation*}\begin{split}
&\Infl = \Set{(t,c) \in [0,1]\times\CC_f \st 
   \LiftSeg{\widehat{0},\fhat(re^{2\pi it})} \cap \Vor{\fhat(c)} \ne
   \emptyset}\\
\Infl_t = &\Set{c \in \CC_f \st (t,c) \in \Infl} 
  \qquad\qquad\qquad
\Infl_c = \Set{t \in [0,1] \st (t,c) \in \Infl}
\end{split}\end{equation*}
\end{defn}

Notice that, for $z = re^{2\pi it}$ fixed, we have  $c\in\Infl_t$
precisely when, for some $y \in \ray{f(z)}$, 
$D_{|f(c)-y|}(y)$ is the largest ball on which $f_z^{-1}$
is defined.  Similarly, for this pair $(t,c)$, we also have 
$t\in\Infl_c$.

\section{The Behavior of $f$ on the Initial Circle} \label{circlebounds}
\bigskip
Consider the function $a_r:[0,1)\to \R$ defined by
$$
a_r(t)=\Arg f(re^{2\pi i t}),
$$
with $r>0$.  We can easily bound the rate of change of $a_r(t)$; while
elementary, these bounds play a crucial role for us.

\begin{lem}[{\bfseries Angular Speed Lemma}]\label{argspeed} 
Let $r>1$. Then for all $t\in [0,1)$, we have
$$
2\pi d\cdot \frac{r}{r+1}\le a_r'(t) \le  2\pi d\cdot \frac{r}{r-1}.
$$
\end{lem}

\begin{proof} Let $z=re^{2\pi i t}$, with $r>1$. Since $|\zeta|\le 1$, we have 
$\frac{\zeta}{z}\in \overline{D_{\frac{1}{r}}(0)}=\Set{w \st |w|\le \frac{1}{r}}$. 
A calculation shows
\begin{equation}\begin{aligned}\label{dar}
a_r'(t)
=\Imm \frac{d}{dt}\log f(re^{2\pi i t})
&=\Imm\left( \frac{d}{dz} \log f(z) \right) \left( re^{2\pi i t} \right)
\cdot 2\pi i \\
&= 2\pi\cdot \Ree\left( \frac{f'(z)}{f(z)}\cdot z \right)
=2\pi\cdot \Ree\sum_{j=1}^d \frac{z}{z-\zeta_j} \\
&=2\pi\cdot \Ree\sum_{j=1}^d \frac{1}{1-\zeta_j/z}.
\end{aligned}\end{equation}

For each root $\zeta_i$, we have
$$ 
\frac{r}{r+1}\le \Ree \frac{1}{1-\zeta_i/z}\le \frac{r}{r-1}.
$$ 
Summing this inequality over the $d$ roots and applying it to
\EqRef{dar} gives the desired result. 
\end{proof}

\begin{rem} The estimates in Lemma~\ref{argspeed} are sharp.
\end{rem}

The following bounds $\alpha(z)$ for points on the initial circle.
This will be of use in proving 
Lemma~\ref{WidthOfStartingPoints}, 
used in selecting starting points to locate all $d$ roots of $f$
in Section~\ref{AllRootsSec}.

\begin{cor}\label{AlphaOutsideDisk}
For $z$ with $|z|=1+1/d$, we have
$$  \left| \frac{f(z)}{f'(z)}\right| < \frac{3}{d}, \quad
    \gamma(z) \le \frac{d(d-1)}{2}, \quad \text{and}\quad
    \alpha(z) < \frac{3}{2}(d-1).$$
\end{cor}

\begin{proof}
Since $r=|z|=1+\frac{1}{d}$, Lemma~\ref{argspeed} gives us $\pi\,d < a_r'$.
From this and the observation that $\Ree(w) \le |w|$, we have
$$
 \pi d <  \left|\frac{f'(z)}{f(z)}\right| \cdot 2\pi (1+\frac{1}{d}) < 
   3\pi \left|\frac{f'(z)}{f(z)}\right|,
\qquad\text{and so}\qquad
  \left|\frac{f(z)}{f'(z)}\right| < \frac{3}{d}.
$$

\medskip
Note that if $\xi_i$ are the $k$ solutions to $f^{(k-1)}(\xi_i) = 0$ (with
multiplicity), then by Lucas' Theorem~\cite{Lu}, we have each $\xi_i$ in the
unit disk and so $|z-\xi_i| \ge 1/d$.  Thus 
$$ \left|\frac{f^{(k)}(z)}{f^{(k-1)}(z)}\right| 
= \left|\sum_{i=1}^k \frac{1}{z-\xi_i} \right| \le d(d-k).$$

Observe that
\begin{align*}
\left| \frac{f^{j}(z)}{j! f'(z)}\right|^{\frac{1}{j-1}} 
&=  \left|\frac{1}{j!} \frac{f''(z)}{f'(z)} \cdot \frac{f'''(z)}{f''(z)} \cdots 
                    \frac{f^{(j)}(z)}{f^{(j-1)}(z)} \right| ^{\frac{1}{j-1}} \\
&\le \left(\frac{1}{j!} d(d-1)\cdot d(d-2) \cdots d(d-j+1) \right)^{\frac{1}{j-1}}
\le \frac{d(d-1)}{2}. 
\end{align*}
Since $\gamma(z)$ is the maximum of the above expression over $j$, we have
$\gamma \le \frac{1}{2}d(d-1)$; combining the two estimates also gives
$\alpha(z) < \frac{3}{2}(d-1)$. 
\end{proof}

\medskip
The corollary below has direct implications for path-following methods
that use a constant ratio step-size (such as \cite{Sm85} or \cite{KS}),
which need a cone of a given angular width about $\ray{w_0}$ containing no
(relevant) critical values in order to set the stepsize that ensures
convergence.  The $\alpha$-step algorithm  considered here adjusts for the
presence of critical values (unless they fall on $\ray{w_0}$) and does not  
need a constant width cone, although a starting value lying in
$\Bad{\theta}$ will have a contribution of at least $\log(1/\theta)$ to the
arithmetic complexity caused by the corresponding critical point $c$.
Recall from Definition~\ref{InflSets} that
$c\in\Infl_t$ means that 
the segment $\LiftSeg{\widehat{0}, \fhat(re^{2\pi it})}\in\SS$ intersects $\Vor{c}$.

\begin{cor}\label{goodangle}
Let $r=1+1/d$, and define
\Oldnote{changed ``all $c$'' to ``some $c$.
   \par Added $c\in\Infl_t$; ref asked if we meant $\fhat$}
$$\Bad{\theta} = \Set{ t\in [0,1) \ST{ 
   \bigl| \Arg \frac{f(re^{2\pi it})} {f(c)} \bigr| < \theta,}
   \ \mbox{for some critical point $c\in\Infl_t$}}.$$
Then
\[ \text{measure}(\Bad{\theta}) \le \frac{2\theta}{\pi}\cdot\frac{d-1}{d} .\]
\end{cor}

\begin{proof}
For fixed $r$, the set $\Bad{\theta}$ consists of the inverse image by
$\fhat$ of $d-1$ arcs of angle $2\theta$ in $\SS$ (one for each critical
point). Each of these will grow by no more than $1/\min a_r'(t)$, so
by Lemma~\ref{argspeed}, when $r=1+1/d$ we have
\Oldnote{ref claims we lost a $2\pi$.  Added comment to clarify. Changed $\phi$
to $\theta$}
$$\text{measure}(\Bad{\theta}) 
                      \le \sum_{c\in\CC_f}\frac{2\theta}{\max a_r'(t)} 
                      \le (d-1)\frac{r+1}{2\pi rd}
                      =  (d-1) \frac {\theta (2d+1)}{\pi d(d+1)} 
                      \le \frac{2\theta (d-1)}{\pi d}.
\qedhere
$$ 
Recall that here we are using the convention that the circle has measure $1$. 
\end{proof}

\begin{rem}\label{speedremark}
Let $\Good{\theta}$ be the complementary notion to $\Bad{\theta}$, that is,
$$\Good{\theta} = \Set{ t\in [0,1) \ST{ 
   \bigl| \Arg \frac{f(re^{2\pi it})} {f(c)} \bigr| \ge \theta,}
   \ \mbox{for all critical points $c\in\Infl_t$}}.$$
%
For each $t \in \Good{\theta}$, 
$f^{-1}_{re^{2\pi it}}\colon \C_{target}\rightarrow \C_{source}$ will
be analytic in a cone 
\Oldnote{added $\in\C_{target}$}
$$ \Set{ w \in \C_{target} \ST{ |\Arg(w) - \Arg(f(re^{2\pi it}))| < \theta}}, $$
and consequently such $t$ correspond to ``good starting points'' for a
path-lifting algorithm: 
\Oldnote{added ``with a fixed step-size''}
in a method with a fixed-ratio stepsize, the convergence is assured, and for
the $\alpha$-step algorithm, convergence is rapid.

This is essentially 
Condition~$\Theta$ of \cite{Sm85} and \cite{SS_Complexity2}, with
$\theta=\pi/12$.  
Both these works use $V_f$ to denote our $\Good{\pi/12}$ (also
taking $r=3/2$), and show in Prop.~2 that $\Good{\pi/12}$ has
measure at least $1/6$.  
Above in Corollary~\ref{goodangle}, we show that the measure
of $\Good{\pi/12}$ is at least $5/6$. 
\end{rem}
\bigskip

Recall 
from Section~\ref{prelim} \Oldnote{added ref. for ``recall''}
that the circle of radius $r$ is denoted by $S_r = \Set{z \st |z|=r}$.

\begin{lem}\label{gamcap}
Let $c$ be a critical point on the boundary of $\Bas(\zeta)$, and let
$\vargamma_c$ be the solution to the Newton flow emanating from $c$
whose interior lies in $\Bas(\zeta)$.  Then if $r> 1$, $\vargamma_c \cap
S_r=\emptyset$.
\end{lem}

\begin{proof}Note that the Newton flow points 
  inward on $S_r$ for $r > 1$, which follows from the
  observation that 
  $$\frac{f(z)}{f'(z)} = \frac{1}{\sum\frac{1}{z-\zeta_i}}.$$
  
  \Oldnote{added discussion of uniqueness of $\gamma_c$.}
The uniqueness of $\gamma_c$ follows from Lemma~\ref{EachBasinCoversPlane}
(which says that $f$ is a biholomorphism from $\Bas(\zeta)$ onto a slit plane)
and the observation that $f$ sends solutions into rays: if there were two
solutions $\gamma_c$ and $\varphi_c$ both emanating from $c$ and lying in
$\Bas(\zeta)$, $f(\gamma_c)$ and $f(\varphi_c)$ would coincide near $0$, and
thus $\gamma_c = \varphi_c$. 

  The transversality and uniqueness facts immediately imply Lemma~\ref{gamcap}.

The transversality of the Newton flow to $S_r$ appears in many places (e.g.,
\cite{STW}), but we include a justification here.
Observe that since $|z|>1$ and $|\zeta_i| \le 1$, the vectors 
$z-\zeta_i$ all lie in a half-plane $\HH$ which does not
include the origin.
Consequently, their inverses and hence their sum  $\sum 1/(z-\zeta_i)$ lie in
a (possibly different) half-plane $\HH'$. Inverting again
gives $f(z)/f'(z) \in \HH$. Since $f(z)/f'(z)$ lies in the
original half-plane $\HH$, it is transverse to $S_r$.
\end{proof}

\Oldnote{rewrote because referee objected to the word ``neck''}
Observe that $\Bas(\zeta) \smallsetminus D_1(0)$ will consist of one or
more connected components.  
The following lemma enables us to estimate the width of these.

\begin{lem}\label{neck} Let $r>1$, $\zeta\in \Roots_f$, and let $\upsilon$ be a
  connected component of $S_r \cap \overline{\Bas(\zeta)}$. 
\Oldnote{added where the minimum is taken.}
Then
$$
\length(\upsilon)\cdot \min~a_r'(t) \le 2 \pi r,
$$
where the minimum is taken over points with $r e^{2\pi it} \in \upsilon$.
\end{lem}

\begin{proof} Let $B\subset \overline{\Bas(\zeta)}$ be a boundary
  component of $\Bas(\zeta)$ 
  which does not intersect $\upsilon$, and let $c$ be a critical point of
  $f$ contained in $B$.  Let $\vargamma_c$ be the orbit of the Newton flow
  which begins at~$c$ and ends at the root $\zeta$; then $\vargamma_c
  \smallsetminus \Set{c}$ will
  be contained in $\Bas(\zeta)$
since $f(\vargamma_c)\in\C_{target}$ is the segment $(0,f(c))$.
  
%

Observe that $f(\vargamma_c\cup B)$ is exactly the ray through $f(c)$.  From the
definition of $\upsilon$ and 
Lemma~\ref{gamcap} we get $\inter(\upsilon)\cap(B\cup 
\vargamma_c)=\emptyset$. Hence, 
$$
\Arg(f(\inter(\upsilon)))\cap \Arg(f(c))=\emptyset,
$$
that is, the image of $\upsilon$ cannot make more than a full turn in
the target space. 
The lemma follows.
\end{proof}

The following corollary follows immediately from the proof.

\begin{cor} \label{neckcor}
Let $z_1$ and $z_2$ satisfy $|z_1|=|z_2|=r$ with $r\ge1$, and suppose
also that they lie in the same
connected component of $S_r \cap \Bas(\zeta)$.  Then there is a well-defined
branch of the argument $\Arg$ which is continuous on 
$S_r \cap \Bas(\zeta)$ and such that
$$|\Arg f(z_1) - \Arg f(z_2)| \le 2\pi.$$
\end{cor}

\medskip
In the sequel we will consider integrals over the circle
$S_r=\Set{z\in \C \st |z|=r}$, 
which, for all  $r>0$, carries Lebesgue measure with unit mass. 

We require the following lemma and its corollary in the proofs of
Lemma~\ref{derivProductBound} and Lemma~\ref{logw0wend}.

\begin{lem}\label{intlog|f|} Let $r>0$ and $|\zeta|<r$ then
$$
\int_0^1 \log |re^{2\pi i t}-\zeta|dt= \log r.
$$
\end{lem}

\begin{proof} Define
$$
\begin{aligned}
S(\zeta)=\int_0^1 \log |re^{2\pi i t}-\zeta|dt 
&= \int_{S_r} \Ree (\log (z-\zeta))\cdot \frac{1}{2\pi i}\frac{dz}{z}\\
&=\Ree \frac{1}{2\pi i}\int_{S_r}\log (z-\zeta)\cdot\frac{dz}{z}.\\
\end{aligned}
$$
Note that
$$\begin{aligned}
\frac{dS}{d\zeta}&=-\Ree  \frac{1}{2\pi i}  \int_{S_r} \frac{1} {z-\zeta}
\frac{dz}{z}\\
&=-\Ree \frac{1}{2\pi i}  \int_{S_r} 
        \left(\frac{1/\zeta} {z-\zeta}-\frac{1/\zeta}{z}\right)dz=0.\\
\end{aligned}
$$
Hence,
\[ S(\zeta)=S(0)=\log r. 
\qedhere \]
\end{proof}

The following corollary is needed in the proof of
Lemma~\ref{logw0wend}, but is also interesting in its own right.
 
\begin{cor}\label{Key1} Let $f(z) = \prod_{j=1}^d(z-\zeta_j)$, with
  $|\zeta_j|<r$. Then
$$
\int_0^1 \log |f(re^{2\pi i t})|dt=d \log r.
$$
\end{cor}
\begin{rem}
Notice that if $r=1+1/d$, we have $d\log r < 1$.
\end{rem}

\begin{proof} 
\[
\int_0^1 \log |f(re^{2\pi i t})|dt 
 =\int_0^1 \log \left|\prod_{j=1}^d (re^{2\pi i  t}-\zeta_j)\right|dt
 =\sum_{j=1}^d \int_0^1 \log |re^{2\pi i t}-\zeta_j|dt
 =d\log r,
\]
where the last equality follows from Lemma~\ref{intlog|f|}.
\end{proof}

\begin{problem}
The previous corollary shows that the average value of $\log |f(z)|$
on $S_r$ is $d\log r$.  Is there a constant $c_r$ independent of $d$ so
that
$$\text{measure}\Set{ t \st \log|f(re^{2\pi i t})| < d\log r} > c_r ?$$
\end{problem}

\bigskip

We now establish a  lower bound on $|w_0|=|f(z_0)|$ for starting
points $z_0$ on the circle $S_r$ with $r>1$.  
We shall use this in Lemma~\ref{lowerwN} to give a lower bound on the size
of our final point $w_N$.  The existence of such a bound should be expected,
since $z_0$ is taken outside the disk containing all the roots; we need this
result in the proof of Theorem~\ref{AverageCostBound} to handle the case
where $z_0$ is already an approximate zero of $f$.

\begin{prop}\label{minabsoff} Let $z\in \Bas(\zeta)$ with $|z|=r>1$.
Then
$$
|f(z)|\ge s_r \cdot \rho_\zeta,
$$
where $\rho_\zeta$ is the radius of convergence of the branch of $f^{-1}$
taking $0$ to $\zeta$, and $s_r<1$.

\medskip\noindent
If $r > 1 + \frac{2\pi}{d}$, $s_r = \frac{1}{4}$.  Otherwise, for
$r=1+\frac{C}{d}$, $s_r$ is the smallest positive solution of 
$$C=8\pi\frac{s}{(1-s)^2}.$$
\end{prop}
\nobreak
\begin{rem}\label{minabsoffbound}
For $0 < C \le 2\pi$, we have $0 < s_r \le \alphaZero$.  For $C=1$, we have
$s_r \approx 0.0369 > \frac{1}{28}$.
\end{rem}

\goodbreak
\begin{proof}
Without loss of generality, we may assume that $\zeta$ is a non-negative
real number.  
Define $l$ to be the radius of the largest disk centered at
$\zeta$  which is mapped univalently into $D_{\rho_\zeta}(0)$, that is, 
$$
D_l(\zeta)\subset f_\zeta^{-1}(D_{\rho_\zeta}(0)).
$$
Observe that these  lie entirely inside $\Bas(\zeta)$.

Applying the Koebe~$\frac14$-Lemma (\EqRef{Koebe14}) to
$f^{-1}_\zeta$, we then obtain
\begin{equation}\label{lowerboundl}
l\ge \frac{1}{|f'(\zeta)|}\cdot \frac{\rho_\zeta}{4}.
\end{equation}

\begin{figure}[htbp]
\begin{center}
\begin{overpic}[width=0.96\hsize,tics=5]{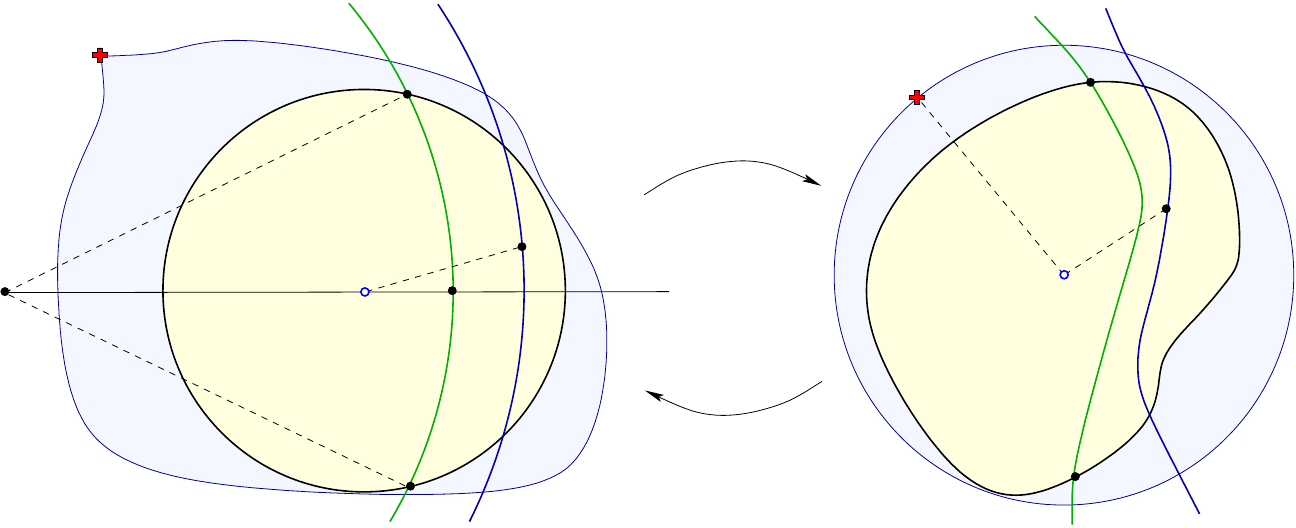}
\put(-0.5,15.5) {\small $0$}
\put(35,16) {\small $1$}
\put(41,22) {\small $z$}
\put(90.5,24) {\small $f(z)$}
\put(6,37) {\small $c$}
\put(67,35) {\small $f(c)$}
\put(72,32) {\small $\rho_\zeta$}
\put(83,23) {\small $s\rho_\zeta$}
\put(5,19) {\tiny $\phi$}
\put(5,16.5) {\tiny $\phi$}
\put(12.5,15) {\small $D_l(\zeta)$}
\put(9,33) {\small $f_\zeta^{-1}(D_{\rho_\zeta}(0))$}
\put(27.5,15.5) {\small $\zeta$}
\put(28,39) {\small $S_1$}
\put(35,39) {\small $S_r$}
\put(75,39) {\small $f(S_1)$}
\put(86,39) {\small $f(S_r)$}
\put(30,31) {\small $A$}
\put(30,4) {\small $\bar{A}$}
\put(81,32) {\small $f(A)$}
\put(78,4) {\small $f(\bar{A})$}
\put(77,17) {\small $f(\zeta)=0$}
\put(56,25) {\small $f$}
\put(55,10) {\small $f^{-1}_\zeta$}
\end{overpic}
\OldnoteC{redid fig~\ref{prflem72} to use \CROSS\ for crit, \CIRC\ for root}
\mycaption{Using the Koebe Lemma to calculate a lower
  bound on $|f(z)|$ for $z$ on $S_r$, in
  Proposition~\ref{minabsoff}.} \label{prflem72} 
\end{center}
\end{figure}

Let $z$ be a point in $\Bas(\zeta)$ with $|z|=r$.  

First consider the case $|z-\zeta| \ge l$.  Here, we must have
$|f(z)|\ge \rho_\zeta/4$. If not, the Koebe $\frac{1}{4}$-Lemma is
violated: by definition of $l$, the map $f$ is univalent on
$D_l(\zeta)$ and so $f(D_l(\zeta))$ contains a disk of radius
$\rho_\zeta/4$ about $0$.
Thus, we need only consider the case when $|z-\zeta|<l$.  

Observe that the function $g(w) = (f^{-1}_\zeta(w) - \zeta)f'(\zeta)$
satisfies the hypotheses of the Koebe Distortion Theorem
(Lemma~\ref{Koebe_Lemma}) on the disk of radius $\rho_\zeta$.
Take $w=f(z)$ to obtain
\begin{equation}\label{UB}
   |z-\zeta| |f'(\zeta)| \le \frac{|f(z)|}{(1-s)^2}
\quad\text{or, equivalently}\quad
   |z-\zeta| \le \frac{1}{|f'(\zeta)|}\cdot  \rho_\zeta \cdot\frac{s}{(1-s)^2},
\end{equation}
where $s=|f(z)|/\rho_\zeta$.

\medskip
We now look for a lower bound on  $ |z-\zeta| $ by estimating
$\frac{|z -\zeta |}{l}$  for   $z \in S_r \cap D_l(\zeta)$. 

\smallskip
Since we have $z \in D_l(\zeta)$ and also $|z|>1$, 
there is a point $A\in S_1 \bigcap D_l(\zeta)$;  let $\phi$
be the angle of the sector connecting $0,A,$ and $1$. See
Figure~\ref{prflem72}.   

Notice that $$ l=\sqrt{ \zeta^2 -2 \zeta \cos(\phi) +1}, 
    \qquad\text{since}\qquad
            (\cos\phi - \zeta)^2 + \sin^2\phi = l^2$$ 
where $(\cos(\phi),\sin(\phi))$  is the coordinate of the point $A$ on
$S_l(\zeta) \cap S_1$.

From Corollary~\ref{neckcor}, we have $|\Arg(f(A))- \Arg(f(\bar{A}))| \le 2\pi$,
and by the Angular Speed Lemma (Lemma~\ref{argspeed}), we have
$$\phi= \Arg(A) \le \frac{ \pi}{d}\cdot \frac{r+1}{r} \le \frac{ 2
\pi}{d}  , \quad \mbox{ for all }  r > 1.$$

Since $r=1+\frac{C}{d}$ and $0< \phi \le \pi$,  we have 
$$\frac{|z-\zeta|}{l}  
\ge \frac{1+\frac{C}{d}-\zeta}{\sqrt{\zeta^2 -2 \zeta \cos(\phi) +1 }}
\ge \frac{1+\frac{C}{d}-\zeta}{\sqrt{\zeta^2 -2 \zeta\cos(\frac{2\pi}{d}) +1}}.
$$
Since we are only considering $0<C<2\pi$ and $|\zeta|\le 1$,  
the above expression is minimized when $\zeta=1$.
Hence, we have
$$
\frac{|z-\zeta|}{l}   \ge
\frac{\frac{C}{d} }{\sqrt{ 1 -2  \cos(\frac{2\pi}{d}) +1 }} \ge
\frac{C}{2 \pi}, 
$$ 
for all  $d$.  Using this with \EqRef{lowerboundl}, we obtain
\begin{equation}\label{bound_on_z_zeta}
|z-\zeta| \;\ge\; \frac{C\, l}{2\pi} \;\ge\; 
          \frac{C}{2\pi}\cdot\frac{\rho_\zeta}{4 |f'\zeta)|}
\end{equation}

This, together with the estimate from \EqRef{UB}, 
gives the lower bound on $s$ as the solution to
$$
\frac{C}{2\pi}\cdot\frac{\rho_\zeta}{4 |f'(\zeta)|} \le 
\frac{s}{(1-s)^2}\,\frac{\rho_\zeta}{|f'(\zeta)|} ,
$$
which simplifies as
$$C \le 8\pi \frac{s}{(1-s)^2}.$$
Denote the smaller positive solution of the above by $s_r$.
Since $s$ was defined by $s=|f(z)|/\rho_\zeta$, this gives us 
$|f(z)|\ge s_r\cdot\rho_\zeta$, as desired.
\end{proof}

\section{The Size of the Step}\label{optimize}

Recall that the $\alpha$-step algorithm (see Section~\ref{alg})
generates a sequence of points $z_n$ with 
\[ z_{n+1} = z_n - \frac{f(z_n) - w_{n+1}}{f'(z_n)} ,\]
where the $w_n$ are a sequence of points tending towards $0$ with the
same argument as $w_0=f(z_0)$.  

In this section, for notational convenience we will
sometimes write $f_n$ for $f(z_n)$, $\fhat_n$ for $\fhat(z_n)$, $f'_n$ for
\Oldnote{added comment about sometimes using $f_n$ etc.}
$f'(z_n)$, $\alpha_n$ for $\alpha(z_n)$, and so on. 

We call the distance between  $w_{n+1}$ and  $w_n$ the
\dfn{$n\Th$-jump} and denote it by 
$$
J_n=|w_{n+1}-w_n|=A\cdot \frac{|f(z_n)|}{\alpha(z_n)} .
$$

\smallskip
The coefficient $A$ (and hence $w_{n+1}$) must be chosen so that $f(z_n)$
will lie close enough to $w_n$ to ensure that the algorithm 
efficiently follows the ray $\ray{w_0}$.  In particular, we show in
Proposition~\ref{Jump}
that taking $A=\frac{1}{15}$ gives us $J_n \ge r_n/66$,
where $r_n$ is the radius of convergence of the appropriate branch of
$f^{-1}$ centered at $w_n$.  The proof of this uses induction; the
inductive hypothesis is established in Proposition~\ref{AInductiveHyp}.
\bigskip

If $f$ were linear, the algorithm would follow $w_n$ exactly, and 
$f(z_n)\equiv w_n$. When the degree of $f$ is at least 2, there will be a small
error which we denote by
$$
\delta_n=|f(z_n)-w_n|.
$$

While the algorithm is described in terms of $\C_{source}$ (the $z_n$) and
$\C_{target}$ ( $f(z_n)$ and the $w_n$), it is more straightforward to think
of it in terms of the branched surface $\SS$. 
\medskip

Let $r_n\ge 0$ be maximal such that 
$$
f_{z_0}^{-1}: D_{r_n}(w_n)\to U
$$ 
is univalent, where $U$ is a neighborhood of $z_n$. This is the distance
between $\widehat{w}_n\in\SS$ and the critical value $v \in \VV_f$ for which
$\widehat{w}_n \in \Vor{v}$. 
Also, let
$R_n\ge 0$ be maximal such that  
$$
f^{-1}_{z_0}:D_{R_n}(f_n)\to V
$$ 
is univalent, where $V$ is a neighborhood of $z_n$. Note that $\fhat_n$
could be in $\Vor{v'}$ for a critical value different from that used for
$\widehat{w}_n$; in this case, we still use $R_n = | v' - f_n|$.

\begin{figure}[htbp]
\begin{center}
\begin{overpic}[width=0.75\hsize,tics=5]{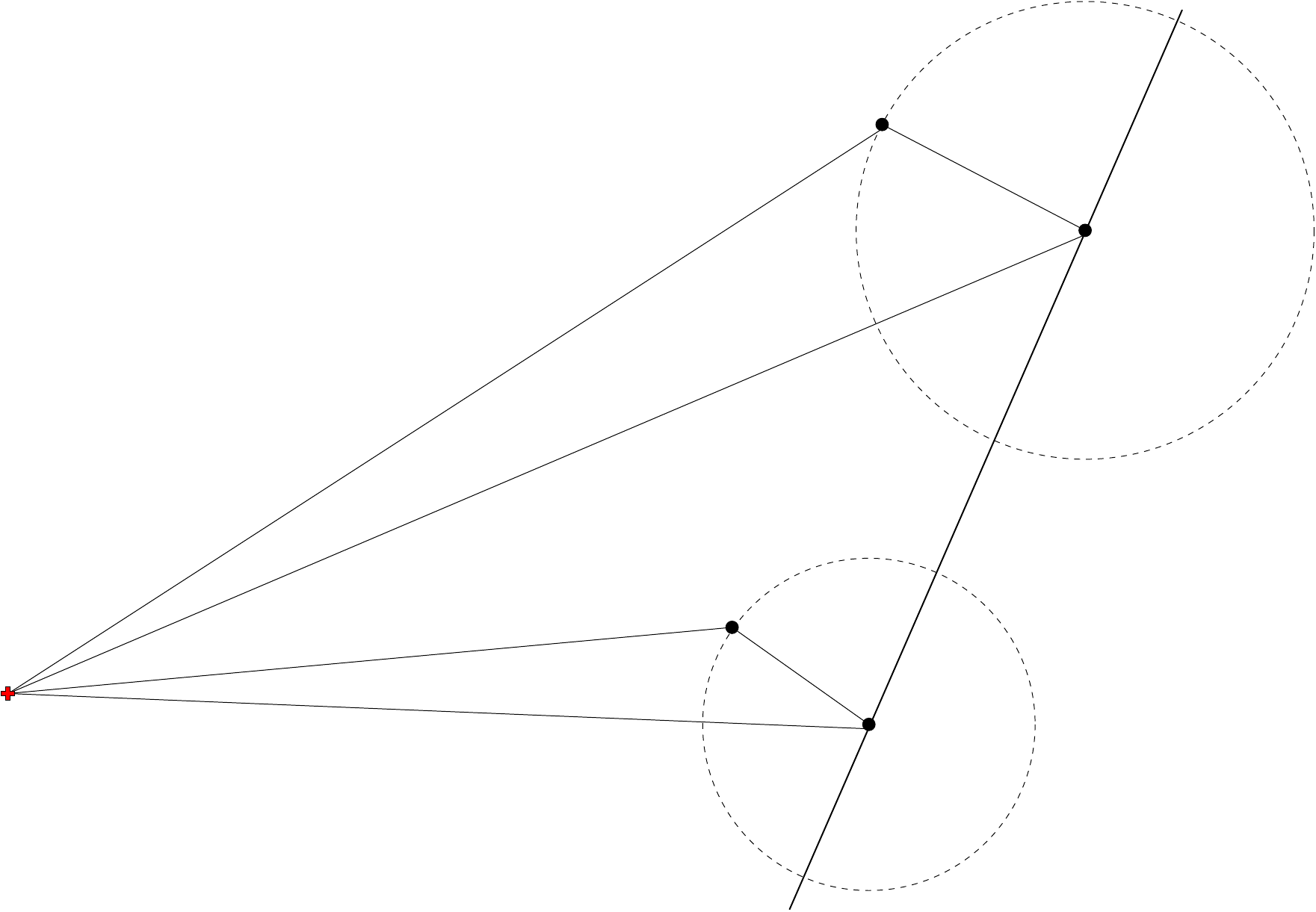}
\put(-2,15)  {\small $v$}
\put(64,61)  {\small $f_n$}
\put(84,51)  {\small $w_n$}
\put(33,40)  {\small $R_n$}
\put(48,39)  {\small $r_n$}
\put(75,57)  {\small $\delta_n$}
\put(52,23.5)  {\small $f_{n+1}$}
\put(67,14)  {\small $w_{n+1}$}
\put(30,21)  {\small $R_{n+1}$}
\put(31,13)  {\small $r_{n+1}$}
\put(60,20)  {\small $\delta_{n+1}$}
\end{overpic}\
\OldnoteC{redid fig~\ref{darR} to use \CROSS\ for crit}
\mycaption{The various notations used througout
  this section, shown in the target space.} \label{darR}
\end{center}
\end{figure}

\bigskip
We introduce the following notation, used throughout this section.
$$
\epsilon_n=z_n - z_{n+1}
\qquad\text{and}\qquad
h_n=(z_n-z_{n+1})\cdot \frac{f'_n}{f_n}=\epsilon_n\cdot \frac{f'_n}{f_n}.
$$
As noted earlier, we use $f_n = f(z_n), f_n' = f'(z_n)$, 
$f_n'' =f''(z_n)$, and $f^{(j)}_n = f^{(j)}(z_n)$ as notation for the derivatives
of $f$ at $z_n$, and use $\alpha_n = \alpha(z_n)$. Let
$\gamma_n=\gamma(z_n)$, where 
\[ \gamma(z)=\max_{j>1}\left|\frac{f^{(j)}(z)}{j!f'(z)}\right|^{\frac{1}{j-1}} \]
as defined in \S\ref{alg}; hence $\alpha_n = \gamma_n |f_n/f_n'|$.

\begin{prop}\label{AInductiveHyp}
Using the preceding notation, suppose we have $A>0$ and $c>0$ given by
\[ \delta_n <  c\cdot\frac{|f_n'|}{\gamma_n} 
   \qquad\text{and}\qquad
   |w_{n+1}-w_n| = A\cdot\frac{|f_n|}{\alpha_n}.
\]
Let $\psi(u)=1 - 4u + 2u^2$.  Then if  $A+c$ satisfies
$(A+c)^2 < c \psi(A+c)^2$, we have
\[ \delta_{n+1} < c\cdot\frac{|f_{n+1}'|}{\gamma_{n+1}}. \]
\end{prop}

\medskip
In order to establish this, we need some preparatory lemmas.

\begin{lem}\label{11} If $|\alpha_n h_n|<1$ then 
$$
\delta_{n+1}=|f_{n+1}-w_{n+1}|
           \le |h_nf_n|\cdot \frac{|\alpha_nh_n|}{|1-\alpha_nh_n|}.
$$
\end{lem}

\begin{proof} Note that since 
$$z_{n+1} = z_n - \frac{f_n-w_{n+1}}{f_n'} ,
\qquad\text{we have}\qquad
w_{n+1} = f_n - (z_n - z_{n+1})f_n' = (1-h_n)f_n.$$
Thus,
\begin{align*}
\delta_{n+1} = |f_{n+1}-(1-h_n)f_n|&=|f(z_n+\epsilon_n)-(1-h_n)f_n|\\
&=\left|f_n+f'_n\epsilon_n+\frac{f''_n}{2!}\epsilon_n^2+\dots-f_n+h_nf_n\right|\\
&=\left|\frac{f''_n}{2!}\epsilon_n^2+\frac{f^{(3)}_n}{3!}\epsilon_n^3+\dots \right|\\
&=|h_nf_n|\cdot 
\left|\frac{f''_n}{2!f'_n}\epsilon_n+\frac{f^{(3)}_n}{3!f'_n}\epsilon_n^2+\dots \right|\\
&\le |h_nf_n|\cdot 
\left|\alpha_n\frac{ f'_n}{f_n}\epsilon_n+(\alpha_n\frac{ f'_n}{f_n}\epsilon_n)^2+\dots \right|\\
&\le |h_nf_n|\cdot 
\left|\alpha_n h_n +(\alpha_n h_n)^2+\dots \right|\\
&\le  |h_nf_n|\cdot \frac{|\alpha_n h_n|}{|1-\alpha_n h_n|}.
\qedhere
\end{align*}
\end{proof}

The proof of the following lemma can be found in \cite{BCSS} 
(Lemma~8.2b and Prop~8.3b).  

\begin{lem} \label{alphaxz}
Let $u_n = \alpha_n h_n$ and $\psi(u) = 1 - 4u + 2u^2$.
Then if $u_n < 1 - 1/\sqrt{2}$, we have
$$ \left| \frac{f_n'}{f_{n+1}'}\right| \le \frac{ (1-u_n)^2}{\psi(u_n)}  
\qquad\mbox{and}\qquad
   \frac{\gamma_{n+1}}{\gamma_{n}} \le \frac{1}{(1-u_n)\psi(u_n)}
$$
\end{lem}

\begin{rem} In \cite{BCSS}, $u_n$ is defined as $(z_n-z_{n+1})\gamma_n$.
  We use $$h_n = \frac{f_n-w_{n+1}}{f_n} = (z_n-z_{n+1})\frac{f_n'}{f_n},$$
  and so our usage and that of \cite{BCSS} agree.  
\end{rem}

\bigskip

We are now ready for the proof of Proposition~\ref{AInductiveHyp}.

\begin{proof}[Proof of Proposition~\ref{AInductiveHyp}]

First, observe that if $A$ and $c$ satisfy 
\begin{equation}\label{AcBasicHyp}
 \delta_n <  c\cdot\frac{|f_n'|}{\gamma_n} 
   \qquad\text{and}\qquad
   |w_{n+1}-w_n| = A\cdot\frac{|f_n|}{\alpha_n}.
\end{equation}
we have $|\alpha_n h_n| \le A+c$.

\begin{align}\begin{split}\label{alphah}
|h_n f_n| &= |f_n-w_{n+1}|\\
&\le |w_n-w_{n+1}|+|f_n-w_n|\\
&\le J_n +\delta_n\\
&\le A \cdot \frac{|f_n|}{\alpha_n}+ c \cdot\frac{|f_n|}{\alpha_n}
 = (A+c)\cdot \frac{|f_n|}{\alpha_n} .
\end{split}\end{align}

We impose the further condition 
%
$$
   A+c < 1 - \frac{1}{\sqrt{2}}
$$
%
which allows us to apply  Lemma~\ref{alphaxz}; this also ensures that the
hypothesis of Lemma~\ref{11} is satisfied.

\medskip

Since $\alpha_n = \gamma_n\cdot|f_n'/f_n|$, by Equation~\ref{alphah}
we have 
$$ |h_nf_n| \le (A + c)\frac{|f_n'|}{\gamma_n}.$$
In Lemma~\ref{11}, we obtained
$$
\delta_{n+1} \le
 \left| h_n f_n \frac{\alpha_n h_n}{1-\alpha_n h_n} \right| \le 
 (A+c)\frac{|f_n'|}{\gamma_n}\cdot \frac{\alpha_n h_n}{1-\alpha_n h_n}.
$$
Thus, it is sufficient to impose the condition
$$
(A+c)\frac{|f_n'|}{\gamma_n}\cdot \frac{\alpha_n h_n}{1-\alpha_n h_n} \le 
c\cdot \frac{|f_{n+1}'|}{\gamma_{n+1}},
$$
or equivalently,
$$
(A+c)\cdot\frac{\gamma_{n+1}}{\gamma_n}\cdot \frac{|f_n'|}{|f_{n+1}'|}\cdot
\frac{1}{c}\cdot \frac{\alpha_n h_n}{1-\alpha_n h_n}<1.
$$
From Lemma~\ref{alphaxz}, after simplification we obtain
$$
  (A+c)\frac{\alpha_n h_n}{\psi(\alpha_n h_n)^2}\cdot\frac{1}{c} < 1.
$$

Since $\alpha_n h_n \le A +c$ and ${u}/{\psi(u)}$ increases monotonically for 
$u\in[0, 1-1/\sqrt{2}]$, we must have 
\begin{equation}\label{Indn+1Fin}
  \frac{(A+c)^2}{\psi(A+c)^2} \cdot \frac{1}{c} < 1.
\end{equation}

Thus, if $A$ and $c$ satisfy the hypotheses of the proposition, the
conclusion follows.
\end{proof}

\begin{rem}\label{AandCvalues}
To optimize the speed of the algorithm, we need to find the largest
$A>0$ for which there is a $c>0$ such that the pair $(A,c)$ satisfies inequality
\EqRef{Indn+1Fin}. Numerics show that such solutions exist for
\Oldnote{Referee asks if ``$A=1/15, c=0.0158$ is rigorously justified''. 
  I updated to $c=1/74$, which is a slightly better estimate. Constants
  updated and checked.}
$A<0.0703039 < 1/14.22396$;  one can readily check that taking $A=1/15$ and 
$c=1/74$ 
satisfies the conditions.
We will use these values of $A$ and $c$ henceforth.
\end{rem}

\bigskip

In order to prove Proposition~\ref{Jump}, we need the following lemma,
which is essentially Corollary~4.3 of \cite{K2}; 
the lower bound of $\frac14$ follows from the Extended L\"owner's Theorem 
in \cite{Sm81}. 
See also \cite{DKST}, where the same constant is obtained for the inverse of
an analytic map between Hilbert spaces.

\begin{lem}\label{alpharadius}
$$
\frac{1}{4}\cdot R_n\le \frac{|f_n|}{\alpha_n} \le \frac{R_n}{\alphaZero}
$$
\end{lem}

\begin{prop}\label{Jump}
If in the $\alpha$-step algorithm, we choose $w_{n+1}$ along $\ray{w_0}$ so that
\[ J_n = |w_n - w_{n+1}| = \frac{1}{15}\cdot \frac{|f_n|}{\alpha_n}, \]
we have $J_n\ge \frac{1}{66} \cdot r_n$ for all $n$.
\end{prop}

\begin{proof}
First, observe that since $w_0=f_0$, we have $\delta_0 = 0$.

Applying Proposition~\ref{AInductiveHyp} with $A=1/15$ and
$c=1/74$ \Oldnote{changed $c$ to 1/74}
then gives us
\begin{equation}\label{deltabound}
\delta_n\le  \frac{1}{74}\cdot\left|\frac{f_n}{\alpha_n}\right|
\end{equation} 
for all $n\ge 0$.

From Lemma~\ref{alpharadius}, we get
$$
J_n =   A\cdot \frac{|f_n|}{\alpha_n} 
    \ge \frac{|f_n|}{15}\cdot \frac{R_n}{4|f_n|}
    =   \frac{1}{60}\cdot \frac{R_n}{r_n} \cdot r_n. 
$$
The radius of convergence at $w_n$ is 
$$
r_n=|w_n-v_n|,
$$
where $v_n$ is the critical value for which $\widehat{w}_n\in\SS$ lies in 
$\Vor{v_n}$. It might be that the radius at $f_n$ is determined by another
critical value, say 
$$
R_n=|f_n-v'_n|.
$$
Let $r'_n=|w_n-v'_n|$. Then we have
$$
r_n \le r'_n  \le |v'_n-f_n|+|f_n-w_n| =   R_n+\delta_n.
$$
In the case when $v_n=v'_n$ we get the same estimate for $r_n$.
Notice, by using \EqRef{deltabound}\ and Lemma~\ref{alpharadius},
$$
r_n \le R_n+\delta_n
\le  R_n+\frac{1}{74}\cdot \frac{|f_n|}{\alpha_n}
\le  R_n+\frac{1/74}{\alphaZero}\cdot R_n
= {\frac{\alphaZero+1/74}{\alphaZero}} \cdot R_n. 
$$        
Consequently, we have
\Oldnote{actually, 65 works instead of 66, but not worth changing}
$$ J_n \ge \frac{\alphaZero}{\alphaZero+1/74}\cdot \frac{r_n}{60} 
   > \frac{r_n}{66},$$
as desired.
\end{proof}

\bigskip
The following corollary tells us how well $f_n$ tracks $w_n$ and
how $w_{n+1}$ relates to $w_n$ as the algorithm progresses.  We use
this below in order to estimate the size of our final guide point $w_N$.

\begin{cor} \label{lowerbw} 
If $\alpha_n > \alphaZero$, then
$$
|f_n| \le \frac{35}{38}\cdot|w_n| 
\qquad\text{and}\qquad
|w_{n+1}|\ge \frac{30}{49}\cdot|w_n|. 
$$
\end{cor}

\begin{proof} Observe,
\Oldnote{replaced constants in cor.\ref{lowerbw}}
\[
|f_n|\le w_n+\delta_n\le |w_n|+\frac{1}{74} \cdot \frac{|f_n|}{\alpha_n}.
\]
Hence,
\[
|f_n|\le \frac{1}{1-\frac{1/74}{\alpha_n}}\,|w_n| =
         \frac{\alpha_n}{\alpha_n-\frac{1}{74}}\,|w_n| \le
         \frac{38}{35}\cdot|w_n|,
\]
where we used $\alpha_n > \alphaZero$ to finish the estimate.

For the second estimate, we have
\begin{align*}
 |w_{n+1}| = |w_n| - \frac{1}{15}\cdot \frac{|f_n|}{\alpha_n} 
 &\ge |w_n| -\frac{1}{15 \alpha_n}\cdot
        \frac{\alpha_n}{\alpha_n-\frac{1}{74}}\,|w_n| \\
 &\ge |w_n| \cdot\left(1-\frac{1}{15}\cdot\frac{1}{\alphaZero - \frac{1}{74}}\right)
  \ge \frac{30}{49}\cdot|w_n|.
\qedhere
\end{align*}
\end{proof}

Using this corollary, we can also obtain a relationship between the guide
point $w_N$ where the algorithm terminates and $\rho_\zeta$, the norm of the
closest critical value to $0$. 
\Oldnote{changed ``recall that'' to ``since\ldots''}
Since the algorithm halts when $w_N$ is an approximate zero for $f$, 
we have 
$\alpha_N \le \alphaZero$ but $\alpha_{N-1} > \alphaZero$.


\begin{lem}\label{lowerwN} For $r\ge 1+\frac1d$ 
$$
|w_N|\ge \frac{1}{40} \cdot \rho_\zeta. 
$$
\end{lem}
 
\begin{proof} From   Proposition~\ref{minabsoff} and
  Remark~\ref{minabsoffbound}, we have 
$$
 |w_0|\ge s_r\cdot \rho_\zeta \ge{\frac{\rho_\zeta}{28}}.
$$

If $w_N = w_0$, the lemma holds trivially.

If $N>0$, then $\alpha_{N-1}\ge \alphaZero$ (and $\alpha_N \le
\alphaZero$). 

\smallskip
From Lemma~\ref{alpharadius}, we get
\Oldnote{I have no idea where 0.032675 came from in the old version.
 I replaced with $\alpha_0/4$ and updated relevant constants, which 
 changes constant from 65 to 40}
\[ |f_{N-1}| \ge \frac14 \cdot  \alpha_{N-1}\cdot R_{N-1}
   \ge \tfrac{\alphaZero}{4} \cdot  R_{N-1} 
   \ge \tfrac{\alphaZero}{4}\cdot \left(\rho_\zeta -|f_{N-1}|\right).
\]
This last inequality follows from the triangle inequality: if $v$ is the
critical value with $|v|=\rho_\zeta$, then $0$, $v$, and $f_{N-1}$ form
a triangle with side lengths $\rho_\zeta$, $R_{N-1}$, and
$|f_{N-1}|$. Rewriting the above yields  
\begin{equation}\label{fN_vs_rho}
|f_{N-1}|\ge \tfrac{\alphaZero}{4+\alphaZero}  \cdot \rho_\zeta.
\end{equation}
We now apply Corollary~\ref{lowerbw} to obtain
\begin{equation} \label{wN_estimate}
|w_N| \ge \frac{30}{49} \cdot |w_{N-1}| 
      \ge \frac{30}{49} \cdot \frac{f_{N-1}}{35/38} .
\end{equation}
Combining equations \EqRef{fN_vs_rho}\ and \EqRef{wN_estimate} gives
$$ |w_N| \ge \tfrac{30\cdot38\cdot(\alphaZero)}{38\cdot49\cdot(4+\alphaZero)}\cdot \rho_\zeta 
   > \frac{\rho_\zeta}{40} .$$
\end{proof}

Finally, we give a lemma which allows us to measure the size of an angular
neighborhood about a point $z_0$ on the initial circle for which the
$\alpha$-step algorithm will lift $\ray{w_0}$.  We use this in
Section~\ref{AllRootsSec}.

\begin{lem}\label{WidthOfStartingPoints}
For $|z_0|>1+\frac{1}{d}$, if
$$\frac{\delta_0}{|f_0|} \le \frac{1}{111d}, \qquad\text{then}\quad
  \delta_0 < \frac{1}{74}\frac{|f_0'|}{\gamma_0}
$$
and the hypotheses of Prop.~\ref{AInductiveHyp} are satisfied at $z_0$.
\end{lem}
\Oldnote{Lemma~\ref{WidthOfStartingPoints} is new for finding all roots.}

\begin{proof}
Since $|z_0|\ge 1+\frac{1}{d}$, Corollary~\ref{AlphaOutsideDisk} gives us
$|f'_0/f_0| > d/3$ and  $\gamma_0<d^2/2$.  Hence,
$$
  \frac{\delta_0}{|f_0|} \le \frac{1}{111d} 
= \frac{1}{37 d^2}\cdot\frac{d}{3} 
< \frac{1}{37d^2}\left|\frac{f_0'}{f_0}\right| .
$$

Thus
$$\delta_0 \le \frac{|f_0'|}{37d^2} = \frac{|f_0'|}{74} \cdot \frac{2}{d^2}
 < \frac{1}{74} \frac{|f_0'|}{\gamma_0}.$$
 \end{proof}

\section{The Pointwise Cost}\label{cost}

In this section we will estimate the number $\Nsteps{f}(z_0)$ of iterates
needed to find an approximate zero starting at $z_0$. We need some
preparation to be able to state the estimate.
To simplify notation and without loss of generality, throughout this section
we shall assume that $\ray{w_0}$ lies along the positive real axis;
\Oldnote{added phrases about change of variables and $\SS$, and
      $\Nsteps{f}=\infty$, since ref said ``Are these restrictions?''}
this can be ensured by an appropriate change of variables.
Furthermore, we shall assume that no relevant critical values of $f$ lie 
on~$\ray{w_0}$ (that is, $\ray{\fhat(z_0)} \in \SS$ is disjoint
  from~$\VV_f$); otherwise, $\Nsteps{f}(z_0)$ will be infinite.
\bigskip

\begin{figure}[htbp]
\begin{center}
\begin{overpic}[width=0.8\hsize,tics=5]{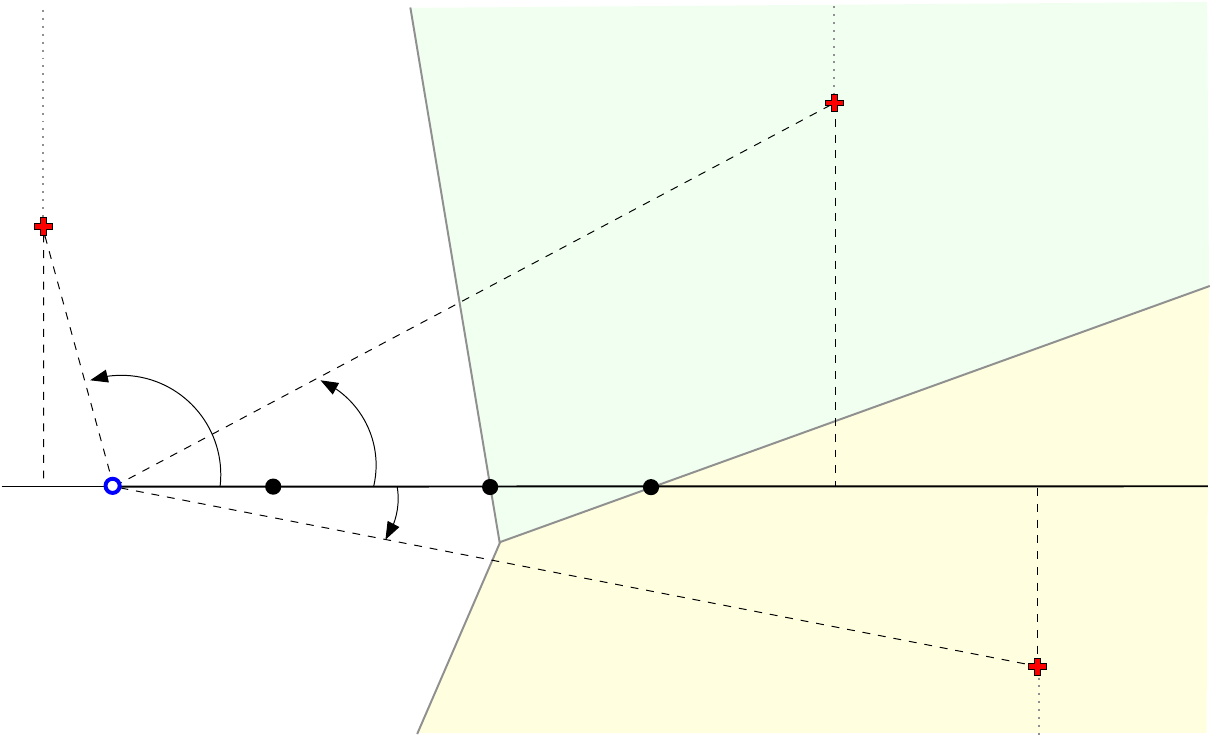}
\put(5,55) {\small $\vor{v_3}$}
\put(36,55){\small $\vor{v_2}$}
\put(37,2) {\small $\vor{v_1}$}
\put(5,42) {\small $v_3$}
\put(70,52) {\small $v_2$}
\put(87,5) {\small $v_1$}
\put(2,18) {\small $p_3$}
\put(67,18) {\small $p_2$}
\put(84,22) {\small $p_1$}
\put(22,23) {\small $q_3$}
\put(40,23) {\small $q_2$}
\put(53,23) {\small $q_1$}
\put(1,30) {\small $x_3$}
\put(66,40) {\small $x_2$}
\put(83,12) {\small $x_1$}
\put(8.5,17) {\small $0$}
\put(101,19) {\small $\ray{w_0}$}
\put(15,28) {\small $\theta_3$}
\put(31,26) {\small $\theta_2$}
\put(33,17) {\small $\theta_1$}
\put(101,19) {\small $\ray{w_0}$}
\end{overpic}
\OldnoteC{redid fig~\ref{defcost} to use \CROSS\ for crit, \CIRC\ for 0}
\mycaption{We divide $\ray{w_0}$ into intervals 
  where it is influenced by each critical value; 
  the various notations used in this section are labeled as in the 
  figure.}\label{defcost}
\end{center}
\end{figure}

As before, let $w_0 = f(z_0)$ and the let the $w_n$ be the guide points along
$\ray{w_0}$ as produced by the algorithm.  
Also let $\widehat{w}_0 = \fhat(z_0)$ and
$\widehat{w}_n$ be the corresponding points in the surface $\SS$, lying along
the ray $\hatray{w_0}$. See Figure~\ref{defcost}.

We divide $\hatray{w_0}$ into subintervals as follows: 
as noted in Proposition~\ref{4pie}, for each $v\in\VV_f$
the intersection of $\hatray{w_0}$ with $\Vor{v}$ will either be an interval or
the empty set.  
Set $\widehat{q}_0=\widehat{w}_0$, and denote the first 
interval by $\LiftSeg{\widehat{q}_0,\widehat{q}_1}$ with
corresponding critical value $v_1$.  In general, set
$$\LiftSeg{\widehat{q}_{j-1},\widehat{q}_j} = \Vor{v_j}\cap\hatray{w_0}.$$
Let $\beta = \beta(z_0)$ denote the total number of such intervals.
Note that for a point $z_0=re^{2\pi it_0}$ on our initial circle, we have
$$ \beta(z_0) = \card \Infl_{t_0},$$
where $\Infl_{t_0}$ is the set of critical points which influence the orbit
of $z_0$, as in Definition~\ref{InflSets}. 

\medskip
So that we may work in the target space $\C$ rather than in the 
surface $\SS$, we make the following observation.
The projection $\pi$ is an isometry in a neighborhood of
$\hatray{w_0}$, since $\VV_f \cap \hatray{w_0} = \emptyset$.
\label{ConstructingU}
\Oldrnote{should we call this about $U$ a remark, to better refer to it in
  Rem.~\ref{znotinbasin}?} 
 We define a set $U(\hatray{w_0}) \subset \SS$ as
$$U(\hatray{w_0}) = 
   \Set{ \widehat{y} \st \LiftSeg{\widehat{y},\widehat{y}_\perp} \ne \emptyset},
$$
where for $y\in\C$, $y_\perp$ denotes the orthogonal
projection of $y$ onto $\ray{w_0}$ (or its extension $\ray{-w_0}$).

That is, for each critical point $c_i$ which influences the orbit of
$w_0$, we remove the ray perpendicular to $\ray{w_0}$ starting at the
critical value $f(c_i)$.  Lifting the result to $\SS$ via the
branch of $\pi^{-1}$ taking $\ray{w_0}$ to $\hatray{w_0}$ yields 
the set $U(\hatray{w_0})$. 

\medskip
Observe that $\pi$ is an isometry on $U(\hatray{w_0})$,
and furthermore, $U(\hatray{w_0})$ contains $\hatray{w_0}$ and a
unique lift of each of the points $f(z_n)$ produced by the algorithm.
Consequently, we have a well-defined correspondence between the target
space $\C$ (minus finitely many rays) and a subset of $\SS$ most
relevant to the $\alpha$-step algorithm starting at $z_0$.  In what
follows, we shall use the notation
  $$\vor{v_i} = \pi(\Vor{v_i}\cap U(\hatray{w_0})),$$ 
and shall slightly abuse notation by using $v_i$ for $f(c_i)$. 

Note that the branch of $f^{-1}$ which takes $w_0$ to $z_0$ is
well-defined throughought all of $\pi(U(\hatray{w_0}))$; in
particular, it coincides with analytic continuation of $f^{-1}$ along
$\ray{w_0}$.  
\medskip

Let $p_j$ be the orthogonal projection of $v_j$ onto the ray
$\ray{w_0}$ (or its extension, $\ray{-w_0}$),  and let $x_j=|v_j-p_j|$.  
See Figure~\ref{defcost}.
Also, let $\theta_j\in (-\pi,\pi]$ be the angle between $v_j$ and the
  ray $\ray{w_0}$; that is,
 $$\theta_j =  \Arg (v_j/w_0).$$ 
Furthermore, use $\beta^+(z_0)$ to denote the number of $\theta_j$ for
which $|\theta_j|\le \pi/2$ (or, equivalently, for which $p_j$ lies on
$\ray{w_0}$). 
\bigskip

With this notation in hand, we can state an upper bound on the cost of
finding an approximate zero starting from a point $z_0$.

\begin{thmA}\refstepcounter{thmcount}\label{PointwiseCostBound}
Let $f\in\PDone$ and let $z_0$ be an initial point for the
$\alpha$-step path-lifting algorithm with $|z_0| > 1$. 
Denote $f(z_0)$ by $w_0$. Then the maximum number of
steps required for the algorithm to produce an approximate zero starting from
$z_0$ is 
\begin{align*}
\Nsteps{f}(z_0) &\le 67 \cdot
 \left( \log\frac{|w_0|}{|w_N|} 
      + \sum_{j=1}^{\beta^+(z_0)} 
        ( 3-2\log|\theta_j| ) 
 \right)
\\
&\le 67 \cdot
 \left( \log\frac{|f(z_0)|}{\rho_\zeta} +\log{40}
      + \sum_{j=1}^{\beta^+(z_0)} 
        ( 3-2\log|\theta_j| )
 \right),
\end{align*}
where $\beta^+(z_0)$ is the number of relevant critical values along
$\ray{w_0}$ with angle $|\theta_j|<\pi/2$, 
and $w_N$ is the final ``guide point'' for the algorithm. 
\end{thmA} 

\begin{rem}
The second inequality follows from the fact that
$\rho_{\zeta}/40 \le |w_N| < \rho_\zeta$, 
as established in Lemma~\ref{lowerwN}.
We shall use this fact in the proving Theorem~\ref{AverageCostBound}.
\end{rem}

\begin{rem} As is shown in  Proposition~\ref{ints} below, for a typical
  starting point,  $\beta^+(z_0) \le 2$  and there are no more than two angles
  $\theta_j$ which are relevant. 
\end{rem}

\begin{rem}\label{algConverges}
In Theorem~\ref{PointwiseCostBound}, the algorithm 
converges to a root $\zeta$ as long as $\theta_j \ne 0$. If $\theta_j=0$,
there is a relevant critical value on $\ray{w_0}$ and the 
algorithm converges to the corresponding critical point; in this case,
$z_0 \not\in \Bas(\zeta)$ for any root $\zeta$ because $z_0$ lies on the 
stable manifold of a critical point.
If $\rho_\zeta=0$, the algorithm will converge to a root $\zeta$ but the
number of steps  $\Nsteps{f}$ will be infinite; in this case $\zeta$ is a
multiple root.  This remark is a restatement of \cite[Thm 5B]{K2} in the
current context.
\end{rem}

In order to establish Theorem~\ref{PointwiseCostBound}, we estimate the
number of steps required to pass each Voronoi domain, and then sum over the
$\beta(z_0)$ domains that $\ray{w_0}$ passes through.  
\medskip

If $w_j$ and $w_k$ are two guide points lying on $\ray{w_0}$ with
$k>j$, we can define the rather trivial function $\Cost{w_j,w_k} =
k-j$.  This measures the number of iterations required by the $\alpha$-step
algorithm beginning at a point $z_j$ near $f_{z_0}^{-1}(w_j)$ to obtain a point
$z_k$ near $f_{z_0}^{-1}(w_k)$.  We extend this function to all pairs of points
$y_1$ and $y_2$ lying on $\ray{w_0}$ by linear interpolation.
It is our goal in this section to estimate $N=\Cost{w_0, w_N}$
where $w_N$ corresponds to an approximate zero of $f$. 

\medskip
Rather than count the number of steps directly (which is possible, but
tedious), instead we follow a suggestion of Mike Shub and integrate
the reciprocal of the stepsize along $\ray{w_0}$.

\begin{lem}\label{costIntegral}
Let $y_1$ and $y_2$ be two points of $\ray{w_0}$.  Then 
$$\Cost{y_1,y_2} \le 67 \int_{y_2}^{y_1}{ \frac{dy}{r_y} },$$
where $r_y = |y-v|$ for each $y\in\vor{v}\cap\ray{w_0}$.
\end{lem}

\begin{proof}
Recall that in Section~\ref{optimize}, we used
$J_n$ to denote the $n\Th$ jump, that is, $J_n = |w_n - w_{n+1}|$
where $w_n$ is a guide point for the algorithm. 
Set $J(w_n)=J_n$, and extend the function $J(y)$ to all of $\ray{w_0}$ by
linear interpolation.
Now consider the differential
equation along $\ray{w_0}$ given by 
\begin{equation}\label{odeJ}
 \frac{dy}{dt} = -J(y) \qquad y(0)=w_0.
\end{equation}
Since $J(y)$ is Lipschitz, \EqRef{odeJ}\ has a unique solution.
Observe that the points $w_n$ are exactly the values given by using Euler's
method with stepsize~1 to solve \EqRef{odeJ} numerically.

Now consider instead the differential equation given by
\begin{equation}\label{odeR}
 \frac{dy}{dt} = -\frac{r_y}{67} \qquad y(0)=w_0.
\end{equation}
We wish to compare the solution of \EqRef{odeR}\ to the Euler method
for \EqRef{odeJ}.  We will show that for every $y$ in any interval 
$[w_{n+1},w_n]$, we have $r_y/67 \le J(y)$.
Consequently, if $\varphi(t)$ is the solution to \EqRef{odeR}\ and
$\varphi(t_1)=y_1$, $\varphi(t_2)=y_2$, then we will have 
$t_2 - t_1 \ge \Cost{y_1, y_2}$.

\medskip

To see that $r_y/67 \le J_y$ for all $y\in [w_{n+1}, w_n]$, we must examine
a few cases.  First, note that if $y\in\vor{v_i}$, we
have 
$$r_y^2 = (y-p_i)^2 + x_i^2.$$
Also, recall that by virtue of Prop.~\ref{Jump}, we have 
$J(w_n) \ge r_{w_n}/66$.

First consider the case where the interval $[w_{n+1}, w_{n}]$ lies entirely
in $\vor{v_i}$.  If $w_{n+1}\ge p_i$, then since $r_y$ is decreasing on the
interval $[p_i, w_n]$, we have $J(y) \ge r_y/66$.  
If $p_i \ge w_{n+1}$, $r_y$ will be nondecreasing.  However, we can apply the
triangle inequality (recalling that $J(w_n)=w_n-w_{n+1}$) to see that 
  $$r_y \le J(w_n) + r_{w_n} \le J(w_n) + 66 J(w_n),$$
and so $J(w_n) \ge r_y/67$ for all $y$ in the interval. 

In the case where the interval intersects more than one Voronoi region, we
proceed as follows.  First, observe that for all $y \in [q_i, w_n]$, we have
already established that $J(y) \ge r_y/67$ holds (where $q_i$ is the
smallest point of $[w_{n+1}, w_{n}]\cap\vor{v_i}$).  Since 
$|v_i - q_i| = |q_i - v_{i+1}|$, we have $J(q_i) \ge r_{q_i}/67$, and
we continue as above.
\medskip

Finally, \EqRef{odeR}\ is separable; elementary calculus yields
\begin{equation*}
t(y) = 67 \int_{y}^{w_0} \frac{dy}{r_y}. \qedhere
\end{equation*}
\end{proof}

Let $y$ be a point on $\ray{w_0}$, and let $c$ be a critical point
which influences $w_0$; as before, let $p$ be the orthogonal
projection of $f(c)$ onto $\ray{w_0}$, and let $x$ denote the 
distance between $f(c)$ and $p$.

\begin{figure}[htbp]
\begin{center}
\smallskip
\begin{overpic}[width=.6\hsize,tics=5]{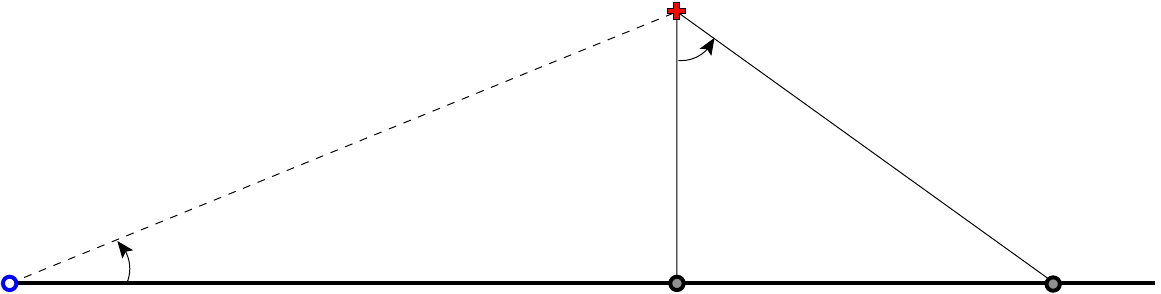}
\put(56,26) {\small $f(c)$}
\put(57,-3) {\small $p$}
\put(90,-3) {\small $y$}
\put(0,-3) {\small  $0$}
\put(12,2) {\small  $\theta_c$}
\put(60,17) {\small  $A_y$}
\put(101,0) {\small $\ray{w_0}$}
\put(55,10) {\small $x$}
\put(73,15) {\small $r_y$}
\end{overpic}
\OldnoteC{redid fig~\ref{triangle} to use \CROSS\ for crit, \CIRC\ for 0}
\mycaption{The quantities $y$, $r_y$, $p$, $x$, $A_y$, and $\theta_c$. }\label{triangle} 
\end{center}
\end{figure}

For each $y$ and a fixed critical point $c$, we also define the angle
$A_y$, which is the angle that the  segment from $y$ to $f(c)$ makes with
the segment between $f(c)$ and $p$.  Notice that $r_y = |f(c)-y|$.
As before, use $\theta_c$ to denote the angle between $f(c)$ and $\ray{w_0}$.
See Figure~\ref{triangle}.

\bigskip

\medskip
We now define the following function, related to $\Cost{y_1, y_2}$:
$$\pounds(y_1,y_2,c) =  
   \log\left( \frac{ (y_1-p) + r_{y_1}}{ (y_2-p) + r_{y_2}} \right).$$ 
By virtue of Lemma~\ref{costIntegral}, if $y_1$ and $y_2$ are both in
$\vor{f(c)}$, we have
\begin{equation}\label{CostInPounds}
\Cost{y_1,y_2} \le 67 \int_{y_2}^{y_1} \frac{dy}{r_y} 
                 = 67\, \pounds(y_1,y_2,c). 
\end{equation}
However, $\pounds$ will still be useful even when one or both of its first two
arguments are not in $\vor{f(c)}$. We establish some bounds
on the value of $\pounds$ in the next few lemmas.

\begin{lem}\label{AyCases}
$$r_y + (y-p) \le 
 \begin{cases} 3(y-p)    & \text{if $A_y > \frac{\pi}{6}$} \\
               x\sqrt{3} & \text{if $A_y \le \frac{\pi}{6}$}
 \end{cases}
$$
\end{lem}
\begin{proof}
Note that $r_y + (y-p) = x (\tan A_y + \sec A_y)$.
If $A_y > \pi/6$, we have $x (\tan A_y + \sec A_y) \le 3 x\,\tan A_y = 3(y-p)$.
When $A_y\le\pi/6$, note that $\tan A_y + \sec A_y$ is increasing in $A_y$;
at $A_y=\pi/6$, $r_y + (y-p) = x\sqrt3$.

We remark that this holds even if $p<0$.
\end{proof}


\begin{lem}\label{y1Toy2}
Let $y_1, y_2 \in \ray{w_0}$ with $y_1 > y_2 \ge 3p >0$.  Then
$$\pounds(y_1,y_2,c) < \log \frac{y_1}{y_2} + \log\frac{9}{4}.$$
\end{lem}
\begin{proof}
We consider two cases: when the angle $A_y$ is large and when it
is small.

If $A_{y_1} \le \pi/6$, since $y_2>p$
$$\pounds(y_1,y_2,c) < \pounds(y_1,p,c) \le \log\frac{x\sqrt{3}}{x} = \log\sqrt{3},$$
where we have used Lemma~\ref{AyCases} in the second inequality.

If $A_{y_1} >\pi/6$, we have (using Lemma~\ref{AyCases} again)
$$\pounds(y_1,y_2,c) \le \log\frac{3(y_1 - p)}{2(y_2 - p)} 
  =  \log\frac{3y_1(1-p/y_1)}{2y_2(1-p/y_2)}.$$
 Since $y_2 \ge 3p$, we have $(1-p/y_1)/(1-p/y_2) < 3/2$, and so
$$\pounds(y_1,y_2,c) \le  \log\frac{y_1}{y_2} + \log\frac{9}{4}.
$$

Since $\sqrt{3} < 9/4$, the above bound holds in either case.
\end{proof}

\begin{lem}\label{3pToZero}
If $p>0$, 
$$
\pounds(3p,0,c) \le  
      \log\frac{4+\tan|\theta_c|}{\sec|\theta_c| -1}. 
$$
\end{lem}
We note that since $p>0$, we have $-\pi/2 < \theta_c < \pi/2$.
Consequently, $\frac{4+\tan|\theta_c|}{\sec|\theta_c| -1} > 1$.

\begin{proof}
We have
$$\pounds(3p,0,c) = \log\frac{(3p - p)+r_{3p}}{r_0-p} \le
      \log\frac{2p + (2p+p\tan|\theta_c|)}{p\sec|\theta_c|-p} = 
      \log\frac{4+\tan|\theta_c|}{\sec|\theta_c| -1}
.$$
\end{proof}

Finally, we handle the case where $|\theta_c|\ge\pi/2$.

\begin{lem}\label{pNegative}
If $y_1 > y_2 > 0 \ge p$, \quad
   $\pounds(y_1,y_2,c) \le \log(y_1/y_2)$.
\end{lem}

\begin{proof}
Observe that $r_{y_2} \ge y_2 - p$, since $r_{y_2}$ is the hypotenuse
of the right triangle with a leg of length $y_2 - p$.  Also, by the
triangle inequality, $r_{y_1} - r_{y_2} \le y_1 - y_2$.

Using this, we have
$$\begin{aligned}
     \frac{r_{y_1} + (y_1 - p)}{r_{y_2} + (y_2 - p)} 
&\le \frac{(r_{y_2} + y_1 -y_2) + (y_1 - p)}{2(y_2 - p)} \\
&=   \frac{2 y_1 -p + r_{y_2} - y_2}{2(y_2 - p)} \\
&\le \frac{2(y_1 - p) + r_{y_2} - (y_2-p)}{2(y_2 - p)} \\
&\le \frac{y_1 - p}{y_2 -p} < \frac{y_1}{y_2}.
\end{aligned}$$

Consequently, 
 $\pounds(y_1,y_2,c) = \log\frac{r_{y_1} + (y_1 - p)}{r_{y_2} + (y_2 - p)} 
  < \log (y_1 / y_2)$ as desired.
\end{proof}

\bigskip
The next lemma enables us to complete the proof of
Theorem~\ref{PointwiseCostBound}. 

\begin{lem}\label{costestimate}
Let $z_0$ be an initial point for the $\alpha$-step path lifting algorithm,
with $|z_0| > 1$, let $f\in\PDone$, $w_0=f(z_0)$. Then the maximum number of
steps required for the algorithm to produce an approximate zero starting from
$z_0$ is 
$$
\Nsteps{f}(z_0) \le 67 \cdot
 \left( \log\frac{|w_0|}{|w_N|} 
      + \beta^+(z_0)\log\frac{9}{4} 
      + \sum_{j=1}^{\beta^+(z_0)} \log
            \biggl(\frac{4+\tan|\theta_j|}{\sec|\theta_j| -1}\biggr)
 \right),
$$
where $\beta^+(z_0)$ is the number of relevant critical values along
$\ray{w_0}$ with angle $|\theta_j|<\pi/2$, 
and $w_N$ is the final ``guide-point'' for the algorithm. 
\end{lem} 

\begin{proof}
First, divide $\ray{w_0}$ into segments where it intersects each of the
$\beta(z_0)$ Voronoi regions $\vor{v_j}$; the $j\Th$ segment will be
bounded by points $q_{j-1}$ and $q_j$ (we set $q_0=w_0$, and
$q_{\beta(z_0)}=w_N$). See Figure~\ref{defcost}.  

Now, we have
\begin{equation}\label{costOfEachInterval}
 N =  \Cost{w_0,w_N} 
   =  \sum_{j=1}^{\beta(z_0)}\Cost{q_{j-1},q_j}
 \le 67 \sum_{j=1}^{\beta(z_0)}\pounds(q_{j-1},q_j,c_j),
\end{equation}
where the inequality follows from Lemma~\ref{costIntegral} and
\EqRef{CostInPounds}.  Applying Lemmas~\ref{y1Toy2} and \ref{3pToZero}
gives us
$$
 \sum_{j=1}^{\beta^+(z_0)}\pounds(q_{j-1},q_j,c_j)
 \le  \sum_{j=1}^{\beta^+(z_0)} \log^+\left|\frac{q_{j-1}}{q_j^*}\right|
      + \beta^+(z_0)\log\frac{9}{4} 
      + \sum_{j=1}^{\beta^+(z_0)} 
           \log \frac{4+\tan|\theta_j|}{\sec|\theta_j| -1}
$$
where $q_j^* = \max(|q_j|,|3p_j|)$.

Note that since $q_j^*\ge |q_j|$, replacing $q_j^*$ with $q_j$ will
still give us an upper bound; furthermore, since $|q_{j-1}|>|q_j|$,
the logarithm of their ratio is positive. 
Thus, we have
\begin{equation}\label{positiveIntervals}
 \sum_{j=1}^{\beta^+(z_0)}\pounds(q_{j-1},q_j,c_j)
 \le  \sum_{j=1}^{\beta^+(z_0)} \log\left|\frac{q_{j-1}}{q_j}\right|
      + \beta^+(z_0)\log\frac{9}{4} 
      + \sum_{j=1}^{\beta^+(z_0)} 
           \log \frac{4+\tan|\theta_j|}{\sec|\theta_j| -1}.
\end{equation}
Now we apply Lemma~\ref{pNegative} to the remaining intervals (if any).
\begin{equation}\label{negativeIntervals}
 \sum_{j=\beta^+(z_0)+1}^{\beta(z_0)}\pounds(q_{j-1},q_j,c_j)
 \le \sum_{j=\beta^+(z_0)+1}^{\beta^(z_0)} \log\left|\frac{q_{j-1}}{q_j}\right|
\end{equation}

Combining \EqRef{positiveIntervals}\ and \EqRef{negativeIntervals}\ 
with \EqRef{costOfEachInterval}\ and recalling that $q_0=w_0$, 
$q_\beta=w_N$ gives the desired result.
\end{proof}

\begin{proof}[Proof of Theorem~\ref{PointwiseCostBound}]
The proof of the main result of this section now follows immediately from 
Lemma~\ref{costestimate}.
First combine the term $\beta^+(z_0)\log\frac{9}{4}$ with the sum, and
then observe that for $|\theta| < \pi/2$, we have
$$  \log \frac{9(4+\tan|\theta|)}{4(\sec|\theta| -1)} 
    \le \log\frac{1}{\theta^2} + 3.$$
This  can be readily seen via the series expansion, which is
$\log(18) - 2\log(\theta)+\theta/4 + \BigO{\theta^2}$.
\end{proof}

\section{The Average Cost}\label{average}

In this section we shall prove Theorem~\ref{AverageCostBound}, 
which follows from averaging the bound found in Section~\ref{cost} 
over the starting points on the circle of radius $r = 1 + C/d$.

\medskip
Recall from Definition~\ref{InflSets} 
that $\Infl$ is the set of pairs $(t,c)$ for which the  
critical points $c \in \CC_f$ influence the starting points 
$z_0 = re^{it}$ on the initial circle of radius $r$, 
$\Infl_t$ is the set of critical points which influence a given
$t$, and $\Infl_c$ are the $t\in S_r$ which are influenced by $c$.

\medskip
For each pair in $(t,c)\in\Infl$, we use $\theta = \theta(t,c)$ to denote the
angle between $[0,f(re^{2\pi it})]$ and $[0,f(c)]$, that is
$$\theta(t,c) = \Arg\frac{f(re^{2\pi it})}{f(c)} .$$
In the notation of Section~\ref{cost}, $\theta(t,c_j) = \theta_j$ where
$v_j = \fhat(c_j)$ and $(t,c_j) \in \Infl$.

Note that for each fixed $c$, $\Infl_c$ is a collection of finitely many
intervals: $\Infl_c$ consists of for those~$t$ such that $\hatray{f(re^{it})}$
intersects $\Vor{\fhat(c)}$.

\bigskip

Define for every critical point $c\in \CC_f$ the function $\theta_c:
\Infl_c \rightarrow \R$ by 
$$
 \theta_c(t)=\theta(t,c) =  \Arg\frac{f(re^{2\pi i t})}{f(c)}.
$$

\begin{lem}\label{twotimes} For each $c\in \CC_f$, the map
  $\theta_c$ is at most $(m_c+1)$-to-one. 
\end{lem}

\begin{proof} For every $\theta\in (-\pi,\pi]$ there are at most
$(m_c+1)$ rays $\hatray{}\subset \SS$ for which the angle between
$[0,f(c)]$ and $\pi({\hatray{}})$ is $\theta$ and which also 
intersect $\Vor{\fhat(c)}$.
This is a consequence  of Proposition~\ref{4pie}. 
\end{proof}

\medskip
As an immediate consequence of the Angular Speed Lemma
(Lemma~\ref{argspeed}), we have 
\begin{equation}\label{argspeedEQ}
2\pi d \cdot \frac{r}{r+1}\le \ \frac{d}{dt}\theta_c(t)\ 
      \le  2\pi d\cdot\frac{r}{r-1}. 
\end{equation}
\medskip

\begin{prop}\label{bound} Let $f \in \PDone$ be of degree $d$ and $r>1$.
 Then
$$
\int_0^1{   \!\!\!\!\!\!
          \sum_{\genfrac{}{}{0pt}{}{c\in \Infl_t}{|\theta(t,c)|<\pi/2}} 
            \!\!\!\!\!\!
          \log\frac{4+\tan|\theta(t,c)|}{\sec|\theta(t,c)|-1} \,dt}
\le 3 \cdot\frac{r+1}{ r}.
$$
\end{prop}

\bigskip
\begin{proof}
Througout the proof, let $\psi(\theta)=\frac{4+\tan|\theta|}{\sec|\theta|-1}$.
From Lemma~\ref{twotimes} and \EqRef{argspeedEQ}, we see that for fixed
values of $c$, we have
$$
\int\limits_{t\in\Infl_c\atop|\theta_c(t)|<\pi/2}\!\! {\log\psi(\theta_c(t)) \,dt}  
\le (m_c+1)\int_{-\pi/2}^{\pi/2}{\log\psi(\theta)\,\frac{d\theta}{\theta_c'(t)} }
\le (m_c+1)\frac{r+1}{2\pi r d}\int_{-\pi/2}^{\pi/2}{\log\psi(\theta)\,d\theta} .
$$
Thus
\begin{align*}
\int_0^1 \sum_{c\in\Infl_t\atop|\theta_c(t)|<\pi/2} \log\psi(\theta(t,c))\,dt
&=  \sum_{c\in \CC_f} \int_{t\in\Infl_c\atop|\theta_c(t)|<\pi/2}\!\!\log \psi(\theta(t,c))\,dt \\
&\le \sum_{c\in \CC_f}(m_c+1)\frac{r+1}{2\pi r d}
          \int_{-\pi/2}^{\pi/2} \log\psi(\theta)\,d\theta \\
&\le \frac{2d -2}{2\pi d}\cdot \frac{r+1}{ r} \cdot 9.2901 \\
&< 3 \cdot \frac{r+1}{r}. 
\qedhere
\end{align*}
\end{proof}

Recall 
from Section\ref{cost} \Oldnote{added ref for ``recall''}
that $\beta^+(z)$ denotes the number of
critical points that influence the orbit of $z=re^{2\pi it}$ with the
critical value in the same half-plane, i.e.,  
$$\beta^+(r e^{2\pi it}) = \card\Set{c\in\Infl_t \st 
                  -\pi/2 < \theta(t,c) < \pi/2}.$$
The next proposition bounds the number of such Voronoi domains a
starting point encounters, on average.

\begin{prop}\label{ints}
$$
\int_0^1 \beta^+(r e^{2\pi i t}) dt \le \frac{1+r}{r}.
$$
\end{prop} 

\begin{proof} 

Note that
$$
\int_0^1 \beta^+(re^{2\pi i t}) \,dt
= \int_0^1 \sum_{c\in\Infl_t \atop |\theta_c(t)|<\pi/2} 1 \,dt
= \sum_{c\in\CC_f} \int_{t\in\Infl_c\atop |\theta_c(t)<\pi/2} 1 \,dt.
$$

As in the proof of Proposition~\ref{bound}, we transport the calculation
from the source space to the target space using the bound on
$\theta_c'(t)$ in \EqRef{argspeedEQ}\ and the fact that
for fixed $c$, $\theta_c(t)$ is at most $(m_c+1)$-to-one
(Lemma~\ref{twotimes}).  This gives us
$$\int_0^1 \beta^+(re^{2\pi i t}) \,dt
 \le \sum_{c\in\CC_f}\int_{-\pi/2}^{\pi/2} \frac{d\theta}{\theta_c'(t)}
 \le \sum_{c\in\CC_f} (m_c + 1) \frac{r+1}{2\pi rd} \cdot \pi
 \le 2(d -1) \frac{r+1}{2rd} 
 < \frac{r+1}{r}. 
$$
Above, we used the fact that $\sum\limits_{c\in\CC_f}m_c = d-1$.
\end{proof}

\bigskip

\begin{lem}\label{logw0wend} If $r\ge 1+\frac{1}{d}$
$$
\int_0^1 \log\frac{|w_0|}{|w_N|} \,dt 
\le  d\log r + \log 40 + \frac{1}{d}\cdot \frac{1+r}{r}\cdot K_f.
$$
\end{lem}

\begin{proof} Corollary~\ref{Key1}, Proposition~\ref{minabsoff},  
Lemma~\ref{argspeed}, and Lemma~\ref{lowerwN} are used in the following
calculation. 
\begin{align*}
\int_0^1 \log\frac{|w_0|}{|w_N|} \,dt
&=   \int_0^1 \log |w_0| \,dt - \int_0^1 \log |w_N| \,dt\\
&\le d \log r-\int_0^1 \log \frac{\rho_\zeta}{40} \,dt\\ 
&\le d \log r +\log 40 +\sum_{\zeta\in \Roots_f}
|\log \rho_\zeta|\cdot  \frac{1}{d}\cdot\frac{1+r}{r}\\
&\le d \log r +\log 40 + \frac{1}{d}\cdot\frac{1+r}{r}\cdot
K_f .
\qedhere
\end{align*}
\end{proof} 

\begin{rem} If $r=1+\frac{1}{d}$, then $d\log r < 1$, giving
$\displaystyle{
\int_0^1 \log\frac{|w_0|}{|w_N|} \,dt 
\le  1 + \log 40 + \frac{2 K_f}{d}}
$.
\end{rem}

\bigskip
Now we are ready to provide a proof of the following

\begin{thmB}\refstepcounter{thmcount}\label{AverageCostBound}
Let $f:\C\to \C$ be a monic polynomial with distinct roots $\zeta_i$ in the
unit disk. Let $\Avgsteps{f}$ be the average
number of steps required by the $\alpha$-step algorithm to
locate an approximate zero for $f$.  Then
$$
\Avgsteps{f}\le 67 \Biggl(12.4 +\frac{2 {K_f}}{d} \Biggr).
$$
where the average is taken over starting points on the circle of radius
$1+1/d$ endowed with uniform measure.  
\end{thmB}

\begin{proof} Let $r=1+1/d$.
Lemma~\ref{costestimate}, Proposition~\ref{ints}, 
Lemma~\ref{logw0wend}, and Proposition~\ref{bound} imply

\begin{align*}
\Avgsteps{f} 
  =& \int_0^1 \Nsteps{f}(r e^{2\pi it}) \,dt \\
\le& \int_0^1  67 \cdot \left[ \log\frac{|w_0|}{|w_N|} + 
               \beta^+(re^{2\pi it})\log\frac{9}{4} + 
               \sum_{c\in\Infl_t \atop |\theta(t,c)|<\pi/2} 
                 \log\frac{4+\tan|\theta(t,c)|}{\sec|\theta(t,c)|-1} \
      \right] \,dt \\
\le& 67 \left[
     \left(1 + \log 40 + \frac{2K_f}{d} \right) +
     1.622 +
     6 \;
     \right ] \\ 
\le& 67 \cdot \left[ 12.4 +  \frac{2K_f}{d} \right].
\qedhere
\end{align*}
\end{proof}


\section{The Relation Between Cost and Degree}\label{CostVsDegree}

In the previous section, we showed that the expected number of steps
required for the algorithm to converge to an approximate zero is
bounded by $\Avgsteps{f}$, which depends directly on $K_f/d$.   
\Oldrnote{ref says ``so we need a separate stopping condition for approx
  multiple roots''  What does he want here?}
For every degree $d$, this is neither bounded above nor
below, even if we restrict $f$ to monic polynomials with distinct roots in 
the unit disk.   As noted in Remark~\ref{KfInfinite}, $K_f$ (and hence
$\Avgsteps{f}$) is infinite precisely when $f$ has a multiple zero.  Since
distinct roots of $f\in\PDone$ may be arbitrarily close together, $K_f$ cannot
be bounded above.


\smallskip
We can, however, estimate the average value of $K_f/d$ as $f$ ranges over
$\PDone$ (in fact, its closure $\PDoneBar$). We shall see in this
section that this average value grows  no faster than linearly in $d$,
using the product measure on the distribution of roots on $\PDoneBar$.
\Oldnote{Added ``using the product measure on roots''. 
  ref asks if we can translate to coefficient measure. Solved by making it
  Question at the end of the section}

\bigskip

The value of $\rho_\zeta$ is closely related to the function $\gamma(z)$
mentioned in Section~\ref{alg}.  Indeed, we have the following
relationship, which enables us to bound $\Avgsteps{f}$ and $K_f$ 
from $\gamma(\zeta)$ and $f'(\zeta)$ at each of the roots $\zeta$.

\begin{lem}\label{gammaVsRho}
Let 
$\displaystyle
\gamma(z)=\max_{j>1}\left|\frac{f^{(j)}(z)}{j!f'(z)}\right|^{\frac{1}{j-1}}$ 
and let $\zeta$ be a nondegenerate root of $f$. 
Then
$$(\alphaZero)\,\frac{|f'(\zeta)|}{\gamma(\zeta)} \le \rho_\zeta
 \le 4\,\frac{|f'(\zeta)|}{\gamma(\zeta)}
.$$
\end{lem}
 
\begin{proof}
This follows immediately from \cite[Theorem~4.1]{K2}.
\end{proof}

\medskip
It is not hard to show by induction that 
$$ f^{(j)}(z) = 
 \sum_{k_1}\sum_{k_2\ne k_1}\sum_{k_3\notin\{k_1,k_2\}}\cdots%
 \sum_{k_j\notin\{k_1,k_2,\ldots,k_{j-1}\}}
 \prod_{\ i\notin\{k_1,k_2,\ldots,k_j\}} (z-\zeta_i),
$$
and so
\begin{equation}\label{boundfj}
 f^{(j)}(\zeta_m) = 
 \sum_{k_2\ne m}\sum_{\ k_3\notin\{m,k_2\}}\cdots\sum_{k_j\notin\{m,k_2,\ldots,k_{j-1}\}}
 \prod_{\ i\notin\{m,k_2,\ldots,k_j\}} (\zeta_m-\zeta_i),
\end{equation}
that is, a sum of $\frac{(d-1)!}{(d-j)!}$ terms, each of which is a
product of $d-j$ factors.  Using this observation, we obtain the
following. (Compare to \cite[Prop.~5.1]{DedieuSep}.)

\begin{lem} \label{gammaAsRootDist}
 $\displaystyle    \gamma(\zeta_m) \le 
   \frac{d-1}{2}\,\dfrac{1}{\min_{i\ne m}|\zeta_m - \zeta_i|}$
\end{lem}
\begin{proof}
Using \EqRef{boundfj}\ above and cancelling common factors
between $f'$ and $f^{(j)}$ yields
\begin{align*}
 \left| \frac{f^{(j)}(\zeta_m)}{j! f'(\zeta_m)} \right| &= 
 \left|\frac{1}{j!} \sum_{k_2\ne m}\sum_{\ k_3\notin\{m,k_2\}}\cdots
                    \sum_{k_j\notin\{m,k_2,\ldots,k_{j-1}\}}
 \dfrac{1}{\prod_{i=k_2,\ldots,k_j} (\zeta_m-\zeta_i) }\right|\\
&\le 
 \frac{1}{j!}\sum_{k_2\ne m}\sum_{\ k_3\notin\{m,k_2\}}\cdots
                     \sum_{k_j\notin\{m,k_2,\ldots,k_{j-1}\}}
 \dfrac{1}{ ( \min_{i\ne m}|\zeta_m-\zeta_i| )^{j-1} } \\
&=
 \frac{1}{d}{d \choose j}
    \left[\dfrac{1}{\min_{i\ne m}|\zeta_m-\zeta_i|}\right]^{j-1}. \\
\noalign{Consequently,} \\
\gamma(\zeta_m) =
\max_{j>1}\left|\frac{f^{(j)}(z)}{j!f'(z)}\right|^{\frac{1}{j-1}} 
&\le 
\max_{j>1}\left| \frac{1}{d}{d \choose j} \right|^{\frac{1}{j-1}}
     \dfrac{1}{\min_{i\ne m}|\zeta_m-\zeta_i|} 
\le \frac{d-1}{2} \dfrac{1}{\min_{i\ne m}|\zeta_m-\zeta_i|}.
\qedhere
\end{align*}
\end{proof}

\smallskip
We now turn to estimating the average value of the components which
control $K_f$: the derivative at each root and the minimal inter-root
distance. 
%
Identify a polynomial $f(z)=\prod_{i=1}^d (z-\zeta_i)$ in $\PDone$ with
the $d$-tuple of its roots, and thus we can view its closure
$\PDoneBar$ as the polydisk $\D^d$.  Using Lebesgue measure  
on $\D^d$ gives $\PDoneBar$ a volume of $\pi^d$.

\begin{lem}\label{avgRootDist}
For each $m$, we have 
$$
\int\limits_{(\zeta_1,\ldots,\zeta_d)\in\D^d}
    \!\!\!\!\!\!
      \log\frac{1}{\min_{i\ne m} |\zeta_m - \zeta_i|}
      \ d\zeta_1\,d\zeta_2 \cdots d\zeta_d 
   ~\le~ 2(d-1) \pi^d.
$$
\end{lem}

\begin{proof}
Without loss of generality, we may take $m=1$.  

Let $|\zeta_1|=R_1$, and let $\zeta_k$ be a root for
which $|\zeta_1 - \zeta_k|$ is minimized.  Set  
$\zeta_k - \zeta_1 = r_k e^{i\theta_k}$.
Let $D_{r_k}(\zeta_1)$ be the disk centered at $\zeta_1$ with
radius $r_k$, and let $E_k = \D\smallsetminus D_{r_k}(\zeta_1)$ denote the part
of the unit disk exterior to it. See Figure~\ref{MHdisks}.
There are two possibilities: either $E_k$ is an annulus (which occurs when $R_1
+ r_k < 1$), or $R_1 + r_k \ge 1$ and $E_k$ is a crescent.
Let $s_k$ represent the arc length of the part of boundary of $E_k$
which contains $\zeta_k$.

\begin{figure}[htbp]
\begin{center}
\begin{overpic}[width=.5\hsize,tics=5]{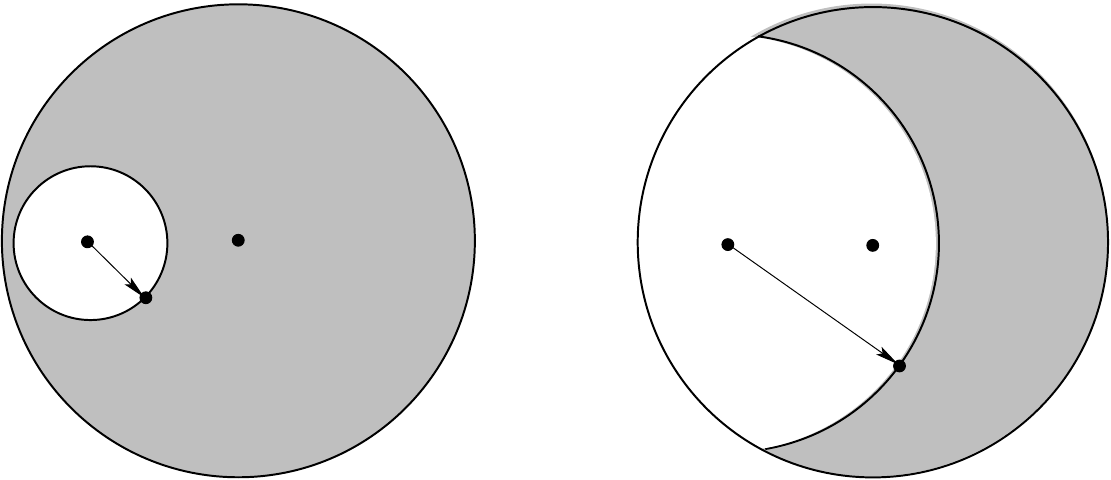}
\put (6,23) {\small $\zeta_1$}
\put (64,23) {\small $\zeta_1$}
\put (20.5,23) {\small $0$}
\put (77.5,23) {\small $0$}
\put (25,35) {\small $E_k$}
\put (85,35) {\small $E_k$}
\put (7,17) {\small $r_k$}
\put (70,14) {\small $r_k$}
\put (82,7) {\small $\zeta_k$}
\put (13,13) {\small $\zeta_k$}
\end{overpic}
\mycaption{The two cases for $E_k$ in Lemma~\ref{avgRootDist}: when
  $r_k+R_1 \le 1$ (left), and when $r_k+R_1 > 1$ (right).  All the
  roots except for $\zeta_1$ lie in the shaded region $E_k$. }\label{MHdisks}
\end{center}
\end{figure}

Observe that for fixed $\zeta_1$, we have 
$(\zeta_2,\ldots,\zeta_d) \in E_k^{d-1}$ 
(with $\zeta_k$ on the interior boundary). 
So we have
$$L_1 =   
\int\limits_{(\zeta_1,\ldots,\zeta_d)\in\D^d}
      \log\frac{1}{\min_{k\ne 1} |\zeta_1 - \zeta_k|}
      \ d\zeta_1\,d\zeta_2 \cdots d\zeta_d 
= \int\limits_{\zeta_1\in\D}\,\int\limits_{(\zeta_2,\ldots,\zeta_d)\in E_k^{d-1}}
    \kern-1em\log\frac{1}{r_k} 
      \ d\zeta_2 \cdots d\zeta_d d\zeta_1 .
$$
The closest root to $\zeta_1$ could be any of remaining $d-1$ roots; 
we shall do the calculation for $\zeta_k$; by symmetry, the remaining
cases will have the same value.

Observe that all roots except $\zeta_1$ lie in in $E_k$. 
The area of $E_k$ is always less than $\pi$ (since it is a subset of
the unit disk), and we always 
have $s_k \le 2\pi r_k$ (since $s_k$ is part of the circumference of a
disk of radius $r_k$.)

If we also write $\zeta_1=R_1e^{i\phi}$ and $\zeta_k-\zeta_1=r_ke^{i\theta_k}$,
and note that integrating $\phi$ and $\theta_k$ give factors of 
$2\pi R_1$ and $s_k$.  Calculating the integral for each $k$ and summing 
gives
$$L_1 \le \pi^{d-2}(d-1)\int_{0}^1
   \int_{0}^{1+R_1} (2\pi R_1) (s_k) \log\frac{1}{r_k} 
   \  dr_k\, dR_1.
$$
Observe that the integrand $\log(1/r_k)$ is positive only for $0<r_k<1$.
Thus, we can give an upper bound on the integral by ignoring the
contribution when $r_k>1$.  

This gives us the following bound on the integral.
\[
L_1 \le 4\pi^d (d-1) \int_0^1 \int_0^1 R_1 r_k\log{\frac{1}{r_k}} \ dr_k dR_1
  = 2(d-1)\pi^d .
\qedhere
\]
\end{proof}

\begin{lem}\label{derivProductBound}
For $f(z)=\prod (z-\zeta_k)$ with $|\zeta_k|\le 1$, we have
 $$\int\limits_{(\zeta_1,\ldots,\zeta_d)\in\D^d}\kern-1.2em
       \log\prod_{m=1}^d\frac{1}{|f'(\zeta_m)|} \ d\zeta_1\cdots d\zeta_d
 ~=~ \frac{d(d-1)}{4} \pi^d .
$$
\end{lem}

\begin{proof}
From \EqRef{boundfj}\ in the case $j=1$, we obtain 
$\displaystyle \prod_{m=1}^d f'(\zeta_m) = \prod_{m=1}^d \prod_{k\ne
  m}(\zeta_m - \zeta_k),$
and so
$$\begin{aligned}
  \int\limits_{(\zeta_1,\ldots,\zeta_d)\in\D^d}\kern-1.2em
  \log\prod_{m=1}^d\frac{1}{|f'(\zeta_m)|} \ d\zeta_1\cdots d\zeta_d
&= -\sum_{m=1}^d\sum_{k\ne m}  \int\limits_{(\zeta_1,\ldots,\zeta_d)\in\D^d}\kern-1.2em
  \log|\zeta_m-\zeta_k| \ d\zeta_1\cdots d\zeta_d \\
&= -\pi^{d-2} \sum_{m=1}^d\sum_{k\ne m} 
   \int\limits_{\zeta_k\in\D}\int\limits_{\zeta_m\in\D} \log|\zeta_m-\zeta_k| 
   \ d\zeta_m d\zeta_k .
\end{aligned}$$

For each of the integrals in the sum, we divide $\D^2$ into two parts:
those where $|\zeta_m| \le |\zeta_k|$ and the complement where
$|\zeta_m|>|\zeta_k|$.   When $|\zeta_m|>|\zeta_k|$, we let 
$\zeta_m = re^{2\pi i t}$ and apply Lemma~\ref{intlog|f|}:
$$\begin{aligned}
   \int\limits_{\zeta_k\in\D}\int\limits_{|\zeta_m|>|\zeta_k|} 
      \log |\zeta_m - \zeta_k|  \ d\zeta_m d\zeta_k
&=  2\pi\kern-.5em\int\limits_{\zeta_k\in\D}\int_{|\zeta_k|}^1\int_0^1 
     \log |re^{2\pi i t}-\zeta_k|   \ r\ dt\,dr\,d\zeta_k &\\
&=  2\pi\kern-.5em\int\limits_{\zeta_k\in\D}\int_{|\zeta_k|}^1
     r\log r \ dr\,d\zeta_k 
&= -\frac{\pi^2}{8}\\
\end{aligned}$$

Similarly, the value of the integral when $|\zeta_m|\le|\zeta_k|$ is also
$-\pi^2/8$. Summing the $d(d-1)$ integrals, each of which contributes
$\pi^d/4$, gives the desired result. 
\end{proof}

\begin{thmC}\refstepcounter{thmcount}\label{LambdaLinearGrowth}
For $f\in\PDoneBar$, let $\Lambda_f$ be the average value of
$\log(1/\rho_\zeta)$, that is, $\Lambda_f = K_f/d$.
Define $\Lambar$ to be the average value of $\Lambda_f$ over
$f\in\PDoneBar$, where we parameterize $\PDoneBar$ by the polydisk of the roots
with Lebesgue measure.  
Then
$$\Lambar < 3d/2.$$
\end{thmC}

\begin{proof}
Applying Lemmas~\ref{gammaVsRho} and \ref{gammaAsRootDist} and using
the fact that $\alphaZero < 1/6$, we have
\begin{align*}
\Lambda_f 
 =   \frac{K_f}{d} 
&=   \frac{1}{d}\sum_{\zeta\in\Roots_f}\log\frac{1}{\rho_\zeta} \\
&\le \frac{1}{d}\sum_{\zeta\in\Roots_f}\log\frac{6\gamma(\zeta)}{|f'(\zeta)|} \\
&\le  \log6 + \frac{1}{d}\sum_{\zeta\in\Roots_f}\log\gamma_\zeta
            + \frac{1}{d}\sum_{\zeta\in\Roots_f}\log\frac{1}{|f'(\zeta)|} \\
&\le  \log6 + \log\frac{d-1}{2} 
            + \frac{1}{d}\sum_{\zeta\in\Roots_f}
                      \log\frac{1}{\min_{\zeta_k\ne \zeta} |\zeta - \zeta_k|}
            + \frac{1}{d}\sum_{\zeta\in\Roots_f}\log\frac{1}{|f'(\zeta)|} .
\end{align*}

Integrating over $f\in\PDoneBar$ and applying Lemma~\ref{avgRootDist} and
Lemma~\ref{derivProductBound}  yields 
$$ \int\limits_{f(z)\in\PDone}\kern-.5em\Lambda_f 
  \le \pi^d \Bigl( \log3 + \log(d-1) +  \frac{2(d-1)}{d} + \frac{d-1}{4} \Bigr)
  < \pi^d\cdot \frac{3d}{2}.
$$
Since the volume of $\PDoneBar$ is $\pi^d$, we obtain
 $\Lambar \le 3d/2$ for all $d$ (and is asymptotic to $d/4$).
\end{proof}

\begin{cor}\label{AvgStepsBound}
For $f\in\PDoneBar$, the average number of steps required to locate an
approximate zero  is $\BigO{d}$.
\end{cor}

\begin{problem}
How does the bound in Theorem~\ref{LambdaLinearGrowth} change if we
average with respect to a measure on the coefficients of~$f$ rather than
uniform measure on the roots of~$f$? 
\end{problem}

\section{How to Find All Roots of a Polynomial}\label{AllRootsSec}

The focus of the paper has been on the question of locating a single
approximate zero for a given polynomial, but these results can easily be used to
locate all $d$ roots of a polynomial $f\in\PDone$.

\smallskip
To do so, we need to locate $d$ initial points, one in
$\Bas(\zeta_j)$ for each root $\zeta_j$.
Then we apply the $\alpha$-step algorithm starting at each of these, and
as long as $f\in\PDone$, the algorithm will produce an approximate zero for
each root.  Our estimates don't rely on roots with special properties (such
as being ``exposed'' as in \cite{Manning}, or having a large sector in the
target space which is free of critical values as in \cite{KS} or
\cite{Sm85}); consequently they apply equally well to each of the roots
$\zeta_j$.

To choose these initial points, we can do the following.
\begin{enumerate}
\item Choose $\left\lceil 111\pi d^2 \right\rceil$ points $y_j$ equally spaced around the circle of radius $1+\frac{1}{d}$.  Let $\tilde{z}_0=y_0$.

\item \label{laststep} 
  Let $k=1$. For each $j>0$, evaluate $f(y_j)$.\\
  If $\Arg f(y_j) \ge \Arg f(y_0)$ but $\Arg f(y_{j-1}) < \Arg f(y_0)$, set
  $\tilde{z}_k=y_j$ and increment $k$. 
\end{enumerate}

At the conclusion of step~(\ref{laststep}), there will be exactly $d$ points
$\tilde{z}_k$ with
$\Arg f(\tilde{z}_k) - \Arg f(y_0) \le \frac{1}{111 d}$. 
This holds as a result of the Angular Speed Lemma (Lemma~\ref{argspeed}) and
the fact that the image of the circle winds exactly $d$ times around the
origin.

Now we use the $d$ points $\tilde{z}_k$ to lift $d$ copies of the same ray 
$\ray{f(y_0)}$, one in each basin, by using a slight modification of the 
$\alpha$-step algorithm from Section~\ref{alg} (Page~\pageref{plm}).
Specifically, we modify Step~0 to set
 $$w_{0,k} = |f(\tilde{z}_k)|\frac{f(y_0)}{|f(y_0)|}, $$
that is, for each $k$ we choose initial target points on the ray
$\ray{f(y_0)}$ with norm $|f(\tilde{z}_k)|$.  Then the $\alpha$-step
algorithm proceeds as usual.

\smallskip
While there could be some $k$ for which $\tilde{z}_k \not\in
\Bas(\zeta_k)$, as a consequence of Lemma~\ref{WidthOfStartingPoints},
each of the points $\tilde{z}_k$ are close enough to some point
$z_{0,k}\in \Bas(\zeta_k)$ (and with $f(z_{0,k}) \in \ray{f(y_0)}$)
so that the $\alpha$-step algorithm will converge to an approximate
zero for the root $\zeta_k$.

The above method for determining the points
$\tilde{z}_k$ requires $\BigO{d^2}$ evaluations of $f$, at an arithmetic
complexity of $\BigO{d^3\log^2 d}$;  the number of steps required
to find all $d$ roots is $\BigO{K_f} = \BigO{\sum \log(1/\rho_f)}$.
Applying Cor.~\ref{AvgStepsBound}, the average complexity to find
approximate zeros for all $d$ roots of $f$ will be $\BigO{d^3\log^2 d}$.

\medskip
\begin{rem}\label{AllRootsRemark}
For $f\in\PDone$, by using the method given above, $d$ approximate zeros
can be found (one for each root $\zeta_j$) in $\BigO{K_f}$ steps of
the $\alpha$-step algorithm.  This has an average arithmetic
complexity of $\BigO{d^3\log^2 d}$.
\end{rem}

\section{Concluding Remarks and Extensions}\label{Remarks}

\begin{rem}\label{HigherNewton}
 Our major goal in this work was to bound the number of iterations of
  the $\alpha$-step algorithm and examine the relationship to the underlying
  geometry of the polynomial, rather than to optimize the arithmetic
  complexity.   Since each step of the algorithm requires computing
  of all of the derivatives of $f$, one could use a higher-order
  method instead of Newton's method (as in \cite{K2},
  \cite{Householder}, \cite{SS_Complexity2})
   in the algorithm without a  significant increase in cost.  In this
   case, we calculate $z_{n+1}$ by a single step of a method using
   higher derivatives of $f$ to approximate the zero of $f(z)-w_{n+1}$
   from $z_n$.
    Use of such a method results in a larger stepsize (and
    consequently fewer steps).  For example, the stepsize is nearly
    doubled by a method using the first three derivatives of~$f$. The
    interested reader should see \cite{K2}, where such methods are examined
    in depth. 
\end{rem}

\begin{rem} \label{NoUsingAlpha}
  Alternatively, the use of $\alpha$ could be curtailed (or even
  entirely removed) by dynamically adjusting the guide points $w_n$ as
  follows.  At each step, set $w_{n+1}$ to be $(1 - h_{n})|f(z_n)| w$.  
  Initially, take $h_{n} = h_0$, but if $f(z_n)$ 
  is not sufficiently close to $w_{n+1}$, divide $h_n$ by 2 and try again 
  until it is. At the next step, start with $h_{n+1} = \min(h_0, 2h_n)$.
  Note that this approach, while similar in spirit, is somewhat
  different from the variable stepsize methods explored in
  \cite{HirschSmale}.
  One can still use $\alpha$ to detect whether an approximate zero has
  been located, or, if evaluating higher derivatives of~$f$ is
  impractical, other methods such as those in \cite{batra_multipoint}
  or \cite{Ocken} can be used.
\end{rem}

\begin{rem} \label{epsilonRootFinding} 
 The $\alpha$-step algorithm could easily be adapted to locate
 $\epsilon$-roots with no significant increase in complexity.  In addition
 to stopping the iteration when an approximate zero is found, the algorithm
 could also stop if $z_n$ is an $\epsilon$-root for a pre-determined
 $\epsilon$.  This can be checked at essentially no cost merely by
 determining if 
  $|f(z_n)/f'(z_n)| < \epsilon/d$ 
(this follows from the well-known fact that there is always a root within 
the disk of radius $d$ times the Newton step at $z$.)
\end{rem}

\begin{rem}\label{Clusters}
\Oldrnote{regarding rem~\ref{Clusters}, ref says ``these ideas could probably
  be used to provide a uniform upper bound on the number of alpha steps
  required (such bounds exist for the Newton dynamics, and in both cases the
  problems are the same, i..e near-multiple roots, so perhaps the treatment
  can be similar?''.  What does he want?} 
 Using some of the ideas in \cite{GLSY}, the results here can be
  extended to deal more directly with multiple roots.
\end{rem}

\begin{rem}
The selection of initial points in Section~\ref{AllRootsSec} can
almost certainly be improved from $\BigO{d^2}$ evaluations of $f$,
most likely to $\BigO{d\log d}$ evaluations. However, this does not
affect the overall complexity of the algorithm.
\end{rem}

\end{document}